\documentclass[12pt]{article}

\usepackage{amsmath}
\usepackage{amssymb}
\usepackage{amsthm}
\usepackage{latexsym}
\usepackage{color}
\usepackage{graphicx}
\usepackage{color}

\DeclareSymbolFont{calletters}{OMS}{cmsy}{m}{n}
\DeclareSymbolFontAlphabet{\mathcal}{calletters}

\def\be{\begin{eqnarray}}
\def\ee{\end{eqnarray}}

\def\b*{\begin{eqnarray*}}
\def\e*{\end{eqnarray*}}

\newtheorem{Theorem}{Theorem}[section]

\newtheorem{Proposition}[Theorem]{Proposition}

\newtheorem{Assumption}[Theorem]{Assumption}

\newtheorem{Lemma}[Theorem]{Lemma}

\newtheorem{Remark}[Theorem]{Remark}
\newtheorem{Example}[Theorem]{Example}

\makeatletter \@addtoreset{equation}{section}

\usepackage{color}

\newcommand{\cblue}{\color{blue}}

\def \E{\mathbb{E}}

\def \L{\mathbb{L}}
\def \P{\mathbb{P}}

\def \R{\mathbb{R}}

\def \M{\mathbb{M}}

\def\Ec{{\cal E}}

\def\Nc{{\cal N}}

\def\Wc{{\cal W}}

\def\Nt{\tilde{N}}

\def \om{\omega}

\def \eps{\varepsilon}
\def \0{\mathbf{0}}

\def \Xh{\widehat{X}}

\def \muh{\hat{\mu}}
\def \sigmah{\hat{\sigma}}
\def \ah{\hat{a}}

\newcommand{\rmi}{{\rm (i)$\>\>$}}

\newcommand{\rmii}{{\rm (ii)$\>\>$}}
\newcommand{\rmiii}{{\rm (iii)$\>\,    \,$}}
\newcommand{\rmiv}{{\rm (iv)$\>\>$}}

\def\x{\times}

\def\1{{\bf 1}}
\def \proof{{\noindent \bf Proof. }}

\def\psih{\widehat{\psi}}

\def\Tt{\tilde{T}}
\def\Wt{\widetilde{W}}

\def\Xt{\widetilde{X}}

\def\Vt{\widetilde{V}}

\def\psit{\widetilde{\psi}}
\def\psib{\overline{\psi}}
\def\Wcb{\overline{\Wc}}
\def\Xb{\overline{X}}

\def\Xbf{\mathbf{X}}

\def\Wch{\widehat{\Wc}}

\def\Thetah{\widehat{\Theta}}
\def\xbf{\mathbf{x}}
\def\mub{\bar{\mu}}

\def\Wct{\widetilde{\Wc}}

\title{Unbiased simulation of stochastic differential equations
\footnote{We are grateful to Emmanuel Gobet, Ahmed Kebaier  and two anonymous referees for valuable remarks and suggestions.
X. Tan and N. Touzi gratefully acknowledge the financial support of the ERC 321111 Rofirm, the ANR Isotace, and the Chairs Financial Risks (Risk Foundation, sponsored by Soci\'et\'e G\'en\'erale) and Finance and Sustainable Development (IEF sponsored by EDF and CA).}}

\author{Pierre Henry-Labord\`ere\thanks{Soci\'et\'e G\'en\'erale, Global Market Quantitative Research,
pierre.henry-labordere@sgcib.com}
        \and Xiaolu Tan\thanks{CEREMADE, University of Paris-Dauphine, PSL Research University, tan@ceremade.dauphine.fr}
        \and Nizar Touzi\thanks{Ecole Polytechnique Paris, Centre de Math\'ematiques Appliqu\'ees,
        nizar.touzi@polytechnique.edu} }

\date{\today}

\begin{document}

\maketitle

\abstract{We propose an unbiased Monte-Carlo estimator for $\E[g(X_{t_1}, \cdots, X_{t_n})]$, where $X$ is a diffusion process defined by a multi-dimensional stochastic differential equation (SDE). The main idea is to start instead from a well-chosen simulatable SDE whose coefficients are updated at independent exponential times. Such a simulatable process can be viewed as a regime-switching SDE, or as a branching diffusion process with one single living particle at all times. In order to compensate for the change of the coefficients of the SDE, our main representation result relies on the automatic differentiation technique induced by Bismu-Elworthy-Li formula from Malliavin calculus, as exploited by Fourni\'e et al. \cite{FLLLT} for the simulation of the Greeks in financial applications.
     In particular, this algorithm can be considered as a variation of the (infinite variance) estimator obtained in Bally and Kohatsu-Higa \cite[Section 6.1]{Bally} as an application of the parametrix method.

    \vspace{5mm}

    \noindent {\bf MSC2010.}    Primary 65C05, 60J60; secondary 60J85, 35K10.

    \vspace{5mm}

    \noindent {\bf Key words.}  Unbiased simulation of SDEs,
    regime switching diffusion,
    linear parabolic PDEs.

}

\section{Introduction}

    Let $d \ge 1$, $T > 0$ and $W$ be a $d$-dimensional Brownian motion,
    $\mu: [0,T] \x \R^d \to \R^d$ and $\sigma: [0,T] \x \R^d \to \M^d$ be the drift and diffusion coefficients,
    where $\M^d$ denotes the collection of all $d \x d$ dimensional matrices.
    Under standard assumptions on these coefficients, we consider the process $X$ defined as the unique strong solution of the multi-dimensional SDE,
    \be \label{eq:SDE}
        X_0 = x_0,
        &\mbox{and}&
        dX_t
        ~=~
        \mu \big(t, X_t \big) ~dt
        ~+~
        \sigma \big(t, X_t\big) ~dW_t,
    \ee
    Our main focus in this paper is on the Monte-Carlo approximation of the expectation
    \be \label{eq:V0_intro}
        V_0 &:=& \E \big[ g\big(X_{t_1}, \cdots, X_{t_n}\big) \big],
    \ee
    for some function $g: \R^{d \x n} \to \R$ and discrete time grid $0 < t_1 < \cdots < t_n = T$. When $n=1$, the analytic formulation of the problem is obtained by the well-known representation $V_0=u(0,X_0)$, where $u$ is the solution of the linear PDE
    \be\label{linearPDE}
    \partial_tu + b(t,x)\cdot Du + \frac12{\rm Tr}\big[\sigma\sigma^\top\!(t,x) D^2u\big]=0,
    &u_T=g,&
    \ee
    In practice, the Monte-Carlo method consists in simulating $N$ independent copies of a discrete-time approximation of $X$, and then estimating $V_0$ by the empirical mean value of the simulations. The corresponding error analysis consists of a statistical error induced by the central limit theorem, and a discretization error which induces a biased Monte Carlo approximation. Under some smoothness conditions, Talay and Tubaro \cite{TalayTubaro} proved that the discretization error for the Euler scheme is controlled by a rate $C \Delta t$, where $\Delta t$ denotes the time step discretization. Since then, many works focused on the analysis of the discretization error under various discretization techniques, see e.g. Kloeden and Platen \cite{KP}, and Graham and Talay \cite{GrahamTalay} for an overview. However, the statistical error estimate $N^{-\frac12}$ is lost in all cases, as its combination with the discretization error leads to an overall error estimate of the order $N^{-\frac12+\eps}$ for some $\eps>0$.

In the context of one-dimensional homogeneous SDEs with constant volatility coefficient, Beskos and Roberts \cite{bes1} developed an exact simulation technique for $X$ by using the Girsanov change measure together with a rejection algorithm, see also Beskos, Papaspiliopoulos and Roberts \cite{bes2}, Jourdain and Sbai \cite{Jourdain}, etc... This technique also applies to more general SDEs by using of the so-called Lamperti transformation which reduces the approximation problem to the unit diffusion case. We also refer to the subsequent active literature of exact simulation of an $\L^\infty-$approximation of $X$, see. e.g Blanchet, Chen and Dong \cite{BlanchetChenDong}.

An alternative  approximation method for $V_0$ was induced by the multilevel Monte-Carlo (MLMC) algorithm introduced by Giles \cite{Giles}, which generalizes the statistical Romberg method of Kebaier \cite{Kebaier}. One of the main advantages of the MLMC algorithm is to control the global error (sum of discretization error and statistical error) with a much better rate w.r.t. the computation complexity. We refer to Giles and Szpruch \cite{GilesSzpruch}, Alaya and Kebaier \cite{AlayaKebaier}, Rhee and Glynn \cite{RheeGlynn} for further developments. In particular, Rhee and Glynn \cite{RheeGlynn} obtained an unbiased simulation method for SDEs by using a random level in the estimator.

More recently,  Bally and Kohatsu-Higa \cite{Bally} provided a probabilistic interpretation of the parametrix method for PDEs. In particular, when $n = 1$, they obtained a representation formula for $V_0$ in form
    \be\label{BKH}
        \E \big[ g\big( \Xh_T \big) \Wc_T \big],
    \ee
    where $\Xh$ is defined by a Euler scheme of $X$ on a random discrete-time grid (the time step follows an independent exponential distribution), and $\Wc_T$ is a corrective weight function depending on $\Xh$.
    The above representation is formally similar to the stochastic finite element method proposed by Bompis and Gobet \cite{BompisGobet}, where one replaces $X$ by its Euler scheme solution in \eqref{eq:V0_intro} and then corrects partially the error by some well-chosen weight functions.
    Notice that in the above representation of  Bally and Kohatsu-Higa \cite{Bally}, the process $\Xh$ can be exactly simulated and hence it may provide an unbiased estimator for $V_0$.
    Nevertheless, the obtained weight function $\Wc_T$ is integrable but has infinite variance, and hence the corresponding Monte-Carlo estimator looses the standard central limit error estimate.

In this paper, we provide a representation of $V_0$ in the spirit of \eqref{BKH}, but with an alternative weight function for the representation. Our results follow from completely different arguments. More importantly, our unbiased approximation of $V_0$ has finite variance, and applies for a large class of SDEs.

Our main idea is to consider the Euler scheme solution $\Xh$ as solution to a regime-switching SDE with some well-chosen coefficients. In order to compensate for the change of the coefficients of the SDE, we introduce some weight functions obtained by the automatic differentiation technique induced by Bismut-Elworthy-Li formula from Malliavin calculus, as exploited by Fourni\'e et al. \cite{FLLLT} for the simulation of the Greeks in financial applications.

The technique introduced in the present paper is inspired by the numerical algorithm introduced in \cite{phl, phl1}, for semilinear PDEs of the form
 \b*
 \partial_t u +\frac12\Delta u + F_0(t,x,u)=0,&u_T=g,
 \e*
for some nonlinearity $F_0$. The main idea in \cite{phl, phl1} is to use an approximation by a branching diffusion representation induced by approximating the nonlinearity $F_0$ by a polynomial in $u$. Namely, given the nature of the linear operator, the representation is obtained by means of a Brownian motion with branching driven by the polynomial approximation of $F_0$.

Loosely speaking, the method developed in the present paper follows by reading the PDE part of \eqref{linearPDE} in the following equivalent form:
 \b*
 \partial_t u +\frac12\Delta u + F_1(t,x,Du,D^2u)
 &=&
 0,
 \e*
where
 \b*
 F_1(t,x,z,\gamma)
 &:=&
 b(t,x)\cdot z+\frac12\mbox{Tr}\big[(\sigma\sigma^\top(t,x)-I)\gamma\big].
 \e*
However, in contrast with the nonlinearity $F_0$, the above function $F_1$ involves the gradient and the Hessian of the solution $u$. Consequently the last PDE can not be handled by the existing literature on branching diffusion representation of PDEs. The automatic differentiation technique introduced in the present paper is an important new idea which allows to convert $Du$ and $D^2u$ in $F_1$ into $u$. Since no powers of $u$ are involved in the equation, this leads to a representation by means of a Brownian motion with exactly one descendent with two different possible types revealed by the weight function corresponding to the order of differentiation.

We believe that the automatic differentiation trick introduced here has very important consequences, beyond the particular application of the present paper. Indeed, in our paper \cite{HenryLabordereOudjaneTanTouziWarin}, we provide a significant extension of the branching diffusion representation to a general class of semilinear PDEs.

    The paper is organized as follows.
    In Section \ref{sec:constant_vol}, we consider the SDE with constant diffusion coefficient, and propose an unbiased estimator for $V_0$ for both Markovian case and path-dependent case.
    Then in Section \ref{sec:general_vol}, we consider the SDE with general diffusion coefficient function, and obtain a similar representation formula for $V_0$, which is integrable but of infinite variance.
    Next, in Section \ref{sec:Numerics}, we provide some numerical examples.
    Finally, we complete some technical proofs in Section \ref{sec:Proofs}.
    In particular, an easy example is studied in Section \ref{subsec:toy_proof} to illustrate the main idea of the technical proofs.

\section{Unbiased simulation of SDE with constant diffusion coefficient}
\label{sec:constant_vol}

	In this section, we will restrict to the constant diffusion coefficient case, 
	and propose an unbiased estimator for $V_0$ having finite variance.

\subsection{The Markovian case}

	Let us start by the Markovian case, where the diffusion process $X$ is defined by
	\be \label{eq:SDE3}
		X_0 = x_0,
		~~~~dX_t
		&=&
		\mu (t, X_t) ~dt ~+~ \sigma_0 ~dW_t,
	\ee
	for some  matrix $\sigma_0 \in \M^d$, and our objective is to compute
	\be \label{eq:def_V0}
		V_0 &=& \E[ g(X_T) ].
	\ee
	for some function $g: \R^d \to \R$. We impose the following conditions on $\mu$ and $\sigma_0$:

	\begin{Assumption} \label{assum:mu_Lip}
    		The diffusion coefficient $\sigma_0$ is non-degenerate,
		the drift function $\mu(t,x)$ is bounded continuous in $(t,x)$,
		uniformly $\frac{1}{2}$-H\"older in $t$ and uniformly Lipschitz in $x$,
		i.e. for some constant $L > 0$,
		\be \label{eq:Lipschitz_coef}
		\Big| \mu(t,x) -\mu(s,y) \Big|
		&\le&
		L \Big( \sqrt{|t-s|} + \big| x-y \big| \Big),
		~~\forall (s,x), (t,y) \in [0,T] \x \R^{d}.~~~~~~
		\ee
	\end{Assumption}

\subsubsection{The unbiased simulation algorithm}

	To introduce our unbiased simulation algorithm, let us first introduce a random discrete time grid.
	Let $\beta > 0$ be a fixed positive constant,
	$(\tau_i)_{i > 0}$ be a sequence of i.i.d. $\Ec(\beta)$-exponential random variables.
	We define
	\be \label{eq:def_T_k}
		T_k ~:=~ \Big( \sum_{i=1}^k \tau_i \Big) \wedge T,~k\ge 0,
		&\mbox{and}&
		N_t ~:=~ \max \big\{ k : T_k < t \big \}.
	\ee
	Then $(N_t)_{0 \le t \le T}$ is a Poisson process with intensity $\beta$ and arrival times $(T_k)_{k> 0}$.
	We denote also $T_0  := 0$ and $\Delta T_{k+1} := T_{k+1} - T_k$.

	Let $W$ be a $d$-dimensional Brownian motion independent of $(\tau_i)_{i > 0}$, we introduce
	\b*
		\Delta W_{T_k} ~:=~  W_{T_k} - W_{T_{k-1}}, ~~~k > 0.
	\e*
	and a process $\Xh$ as the Euler scheme of $X$ on the random discrete grid $(T_k)_{k \ge 0}$, i.e.
	 $\Xh_0 = x_0$ and
	\be \label{eq:REulerScheme}
		\Xh_{T_{k+1}}
		&:=&
		\Xh_{T_k}
		~+~
		\mu\big(T_k, \Xh_{T_k} \big) \Delta T_{k+1}
		~+~
		\sigma_0 \Delta W_{T_{k+1}},
		~~k = 0, 1, \cdots, N_T.
	\ee
	Then our estimator is given by
	\be \label{eq:def_psih1}
		\psih
		~:=~
		e^{\beta T}
		~\Big[ g \Big(\Xh_T \Big)- g\Big( \Xh_{T_{N_T}} \Big) \1_{\{N_T > 0\}}  \Big]
		~\beta^{- N_T}
		~\prod_{k=1}^{N_T}  \Wcb^1_k,~~~
	\ee
	with
	\be \label{eq:def_Wcb}
		\Wcb^1_k
		&:=&
		\frac{ \big( \mu(T_k, \Xh_{T_k}) - \mu(T_{k-1}, \Xh_{T_{k-1}}) \big) \cdot (\sigma_0^{\top})^{-1} \Delta W_{T_{k+1}}}
		{\Delta T_{k+1}}.
	\ee

	\begin{Theorem} \label{theo:renorm_constant_vol}
		Suppose that Assumption \ref{assum:mu_Lip} holds true, and $g$ is Lipschitz.
		Then
		\b*
			\E \big[ \big(\psih \big)^2 \big]
			~<~
			\infty
			&\mbox{and}&
			V_0
			~=~
			\E \big[~ \psih ~\big].
		\e*
	\end{Theorem}
	\proof 
	\rmi We first show that $\E \big[ \big(\psih \big)^2 \big] < \infty$.
    For simplicity, we denote $\Delta \Xh_k := \Xh_{T_k} - \Xh_{T_{k-1}}$ for $k > 0$.
    Let $L_g$ be  the  Lipschitz constant of the function $g$,
    and set $L_0 := \big| \big( \sigma_0 \sigma_0^{\top} \big)^{-1} \big| > 0$ by the non-degeneracy of $\sigma_0$.
    Then using Assumption \ref{assum:mu_Lip},
    it follows by direct computation that
    \b*
        \big| e^{- \beta T} \psih \big|
        \le
        L_g \Big ( |g(x_0)| + \Delta T_1 + |\Delta \Xh_{T_1}| \Big)
        \prod_{k=1}^{N_T} \frac{L (\sqrt{\Delta T_{k+1}} + |\Delta \Xh_{T_{k+1}} | )}{\beta \Delta T_{k+1}}
        \Big| (\sigma_0^{\top})^{-1} \Delta W_{T_{k+1}} \Big|.
    \e*
	Then denoting $\widehat{\E}_{T_k} := \E \big[ \cdot \big|  \Xh_{T_k}, \Delta T_{k+1}\big] $, we have
	\b*
		\widehat{\E}_{T_k} \Big[ \Big| \frac{ \sqrt{\Delta T_{k+1}} + |\Delta \Xh_{{T_k+1}}| }{ \Delta T_{k+1}} (\sigma_0^{\top})^{-1} \Delta W_{T_{k+1}} \Big|^2  \Big] 
		\le
		\E \Big[ \big( 1 + |\mu|_{\infty} \sqrt{T} + \big| \sigma_0 Z \big| \big)^2 \big| (\sigma_0^{\top})^{-1} Z\big| ^2 \Big] ,
	\e*
	where $|\mu|_{\infty} := \sqrt{ \sum_{i = 1}^d |\mu_i|_0^2}$, $|\mu_i|_0 := \sup_{t,x} |\mu_i(t,x)|$,
	and $Z$ is a standard centered normal distribution in $\R^d$.
	This provides 
	\b* 
	&&
	\widehat{\E}_{T_k} \Big[ \Big| \frac{ \sqrt{\Delta T_{k+1}} + |\Delta \Xh_{{T_k+1}}| }{ \Delta T_{k+1}} (\sigma_0^{\top})^{-1} \Delta W_{T_{k+1}} \Big|^2  \Big] \\
	&\le&
	2 \big( 1 + |\mu|_{\infty} \sqrt{T} )^2 ~\E \big[ \big| (\sigma_0^{\top})^{-1} Z \big|^2 \big]
	~+~
	2 \E \big[ \big| \sigma_0 Z  \big|^2 \big| (\sigma_0^{\top})^{-1} Z \big|^2 \big]  \nonumber\\
	&=&
	2  \big( 1 + |\mu|_{\infty} \sqrt{T} )^2 ~\mathrm{Tr} ( (\sigma_0 \sigma_0^{\top})^{-1})
	~+~
	2 \big( 3 d + d(d-1) \big)
	~=:~ \gamma.
    \e*
	We therefore get the following upper bound:
    \be \label{eq:VarianceBound}
        \E \big[ \psih^2 \big]
        ~\le~
        C
        e^{2 \beta T} e^{-\beta T+ \frac{\gamma L^2  T}{ \beta}},
        ~\mbox{where}~
        C := L_g^2 \E \big[ \big( |g(0)| + \Delta T_1 + | \Delta X_1| \big)^2 \big].
    \ee    

	\noindent \rmii The equality $V_0 = \E [ \psih ]$ will be proved in Section \ref{sec:Proofs}, 
	with illustration of the main idea in Section \ref{subsec:toy_proof}.
	\qed

\subsubsection{On the choice of $\beta$}
\label{subsubsec:choice_beta}

	Notice that the random variable $\psih$ in \eqref{eq:def_psih1} can be exactly simulated from a sequence of Gaussian $\Nc(0,1)$ and exponential $\Ec(\beta)$ random variables.
	Then the integrability and representation results in Theorem \ref{theo:renorm_constant_vol}
	induce an unbiased simulation Monte-Carlo method to approximate $V_0$,
	with error induced by the standard central limit theorem.
	
	We next observe that the constant $\beta > 0$ may be chosen so as to minimize the approximation error relative to the computational effort:
	\begin{itemize}
	\item By the central limit theorem, the error induced by the Monte Carlo estimator based on the representation $\hat\psi$ is characterized by the variance of $\psih$. For tractability reasons, we shall instead replace it by the bound \eqref{eq:VarianceBound}.
	\item The computation effort is proportional to the number $N_T$ of arrivals of the Poisson process before the maturity $T$, and is thus given by $C'\E[ N_T] = C'\beta T$.
	\end{itemize}
In view of this, we shall choose $\beta$ by minimizing the ratio of the variance bound \eqref{eq:VarianceBound} to the mean computational effort. This minimization problem is obviously independent of the constants $C,C'$, and reduces to:
	\b*
		\min_{\beta >0 } f(\beta),
		&\mbox{where}&
		f(\beta)
		:=
		\frac{1}{\beta T} \exp \Big( T \big( \beta + \frac{\gamma L^2}{\beta} \big) \Big).
	\e*    
    Direct computation shows that the equation $f'(\beta) = 0$ has a unique solution on $(0, \infty)$ given by
    \b*
        \beta^* &:=& \sqrt{\gamma L^2 + T^2/4} ~+~ \frac{T}{2}.
    \e*
    As $\lim_{\beta \searrow 0} f(\beta) = \lim_{\beta \to \infty} f(\beta) = \infty$, this shows that $\beta^*$ is the minimizer of the above defined criterion, and will be taken as our ``best sub-optimal'' choice of $\beta$ for the unbiased estimator $\psih$.

\subsection{The path-dependent case}

\label{subsec:path_depend}

    In this part, we would like to provide an extension of the above estimator $\psih$ in \eqref{eq:def_psih1}
    to the path-dependent case.
    Let $n > 0$, $0 = t_0 < t_1 < \cdots < t_n = T$, $\sigma_0 \in \M^d$ be a non-degenerate matrix,
    and $\mu: [0,T] \x \R^{d \x n} \to \R^d$ be a continuous function, Lipschitz in the space variable.
    Let $X$ be the unique solution of SDE,
    with initial condition $X_0 = x_0$,
    \be \label{eq:SDE2}
        dX_t
        &=&
        \mu(t, X_{t_1 \wedge t}, \cdots, X_{t_n \wedge t}) ~dt
        ~+~
        \sigma_0 ~dW_t;
    \ee
    and the objective is to compute the value,
    \be \label{eq:def_Vt0}
        \Vt_0
        &:=&
        \E \big[  g \big( X_{t_1}, \cdots, X_{t_n} \big) \big],
    \ee
    for some Lipschitz function $g : \R^{d \x n} \to \R$.

    \begin{Remark} \label{rem:PDE_n}
        It is clear that the value $\Vt_0$ defined above can be characterized by a parabolic PDE system.
        Namely, for every $k = 1, \cdots, n$ and $(x_1, \cdots, x_{k-1}) \in \R^{d \x (k-1)}$,
        we define
        \be \label{eq:def_mu_k}
            \mu_k(t,x) ~:=~ \mu(t,x_1, \cdots, x_{k-1},  x, \cdots, x),
            ~~~\forall (t,x) \in  [t_{k-1}, t_k] \x \R^d.
        \ee
        Suppose that $(u_k)_{k = 1, \cdots, n}$ is a family of functions such that
        $u_k$  is defined on $[t_{k-1}, t_k] \x \R^{d \x k}$ and
        $x \mapsto u_k(t, x_1, \cdots, x_{k-1}, x)$ is a solution (at least in the viscosity sense) of
        \be \label{eq:PDE_n}
            \partial_t u_k ~+~ \frac{1}{2} \sigma_0 \sigma_0^{\top} : D^2 u_k ~+~  \mu_k \cdot D u_k
            &=&
            0,
        \ee
        with terminal conditions
        $$
            u_k(t_k, x_1, \cdots, x_k) = u_{k+1}(t_k, x_1, \cdots, x_k, x_k),
            ~~\mbox{for}~~
            k = 1, \cdots, n-1,
        $$
        and $u_n(t_n, x_1, \cdots, x_n) = g(x_1, \cdots, x_n)$.
        Then we have $\Vt_0 = u_1 (0, x_0)$.


    \end{Remark}

\subsubsection{The algorithm}

    The unbiased simulation algorithm of $\Vt_0$ can be obtained by an iteration of the estimator \eqref{eq:def_psih1}
    on each time interval $[t_k, t_{k+1}]$.
    One should just be careful on the integrability issue.
    Let us first introduce the algorithm.

    Recall that $W$ be a standard $d$-dimensional Brownian motion,
    $(\tau_i)_{i > 0}$ is a sequence of i.i.d. $\Ec(\beta)$-exponential random variables independent of $W$.
    Then $N = (N_s)_{0 \le s \le t}$ and $(T_i)_{i > 0}$ are defined in \eqref{eq:def_T_k}.
    Define further for every $k = 1, \cdots, n$,  $ \Nt^k :=  N_{t_k} -  N_{t_{k-1}}$ the number of jump arrivals on $[t_{k-1}, t_k)$, and $\Tt^k_0 := t_{k-1}$ and $\Tt^k_j := T_{ N_{t_{k-1}} + j} \wedge t_k$,
    \b*
        \Delta \Tt^k_j := \Tt^k_j - \Tt^k_{j-1},
        ~~
       	\Wt^k_j ~:=~ W_{\Tt^k_j},
        &
        \Delta \Wt^k_j := \Wt^k_j - \Wt^k_{j-1},
        &\forall j = 1, \cdots, \Nt^k + 1.
    \e*

    \begin{Example}
        We give below an example for the case $n = 2$.
        In the following example, the number of jump arrivals on $[0, t_1)$ is $ \Nt^1 = 2$,
        that on $[t_1, t_2)$ is $ \Nt^2 = 1$, and total number of jump arrivals is $N_T = 3$.

        For $k =1$, we have $\Tt^1_0 = 0$, $\Tt^1_1 = T_1$, $\Tt^1_2 = T_2$ and $\Tt^1_3 = t_1$;
        $\Wt^1_0 = 0$, $\Wt^1_1 = W_{T_1}$, $\Wt^1_2 = W_{T_2}$ and $\Wt^1_3 = W_{t_1}$.
        For $k =2$, we have $\Tt^2_0 = t_1$, $\Tt^2_1 = T_3$, $\Tt^2_2 = t_2$,
        and $\Wt^2_0 = W_{t_1}$, $\Wt^2_1 = W_{T_3}$ and $\Wt^2_2 = W_{t_2}$.

    \end{Example}

                    \setlength{\unitlength}{0.8cm}
                \begin{picture}(20,1.5)
                        \thicklines

                        \put(0.5,0.5){\vector(1,0){15}}

                        \put(0.92,0){\cblue $0$}
                        \put(1,0.5){\line(0,1){0.15}}

                        \put(3.72,0){\cblue $T_1$}
                        \put(4,0.5){\line(0,1){0.15}}

                        \put(6.72,0){\cblue $T_2$}
                        \put(7,0.5){\line(0,1){0.15}}

                        \put(8.4,0){\cblue $t_1$}
                        \put(8.5, 0.5){\line(0,1){0.15}}

                        \put(10.7 ,0){\cblue $T_3$}
                        \put(11, 0.5){\line(0,1){0.15}}

                        \put(13.82,0){\cblue $t_2$}
                        \put(14, 0.5){\line(0,1){0.15}}

        \end{picture}

    \vspace{6mm}

    We next introduce a process $ \big(\Xt^{k,\xbf}_j \big)$, $\forall j = 0, 1, \cdots, N_k +1$,
    for each $k = 1, \cdots, n$ and initial condition $\xbf = (x_0, x_1, \cdots, x_{k-1}) \in \R^{d \x k}$ by
     $\Xt^{k,\xbf}_0 := x_{k-1}$ and
    \b*
        \Xt^{k,\xbf}_{j+1}
        &:=&
        \Xt^{k,\xbf}_j
        ~+~
        \mu_k\big(\Tt^k_j, \Xt^{k,\xbf}_j \big) \Delta \Tt^k_{j+1}
        ~+~
        \sigma_0 \Delta \Wt^k_{j+1}.
    \e*
    Similarly, for every $j = 1, \cdots, N_k$, we define a automatic differentiation weight, with $\mu_k$ defined by \eqref{eq:def_mu_k},
    \b*
        \Wct^k_j
        &:=&
        \frac{\big(\mu_k \big(\Tt^k_j, \Xt^{k,\xbf}_j \big) - \mu_k \big(\Tt^k_{j-1}, \Xt^{k,\xbf}_{j-1} \big) \big)
        \cdot
        \big(\sigma_0^{\top}\big)^{-1} \Delta \Wt^k_{j+1}}{ \Delta \Tt^k_{j+1}}.
    \e*

    We now introduce the algorithm for the path-dependent case, in a recursive way.
    First, for $\xbf = (x_0, x_1, \cdots, x_n) \in \R^{d \x (n+1)}$, set
	$\psit^{\xbf}_{n+1} := g(x_1, \cdots, x_n)$.
	Next, for $k = 1, \cdots, n$, denote
    $$
        \Xbf^{k,\xbf}
        :=
        (x_0, x_1, \cdots, x_{k-1}, ~\Xt^{k,\xbf}_{\Nt^k+1})
	~~\mbox{and}~~
        \Xbf^{k,\xbf, 0}
        :=
        (x_0, x_1, \cdots, x_{k-1}, ~\Xt^{k,\xbf}_{\Nt^k} \1_{\{\Nt^k > 0 \}}).
    $$
    Then given $\psit^{\cdot}_{k+1}$, we define
    \be \label{eq:def_psib}
        \psit^{\xbf}_k
        &:=&
        e^{\beta(t_k - t_{k-1})}
        \Big(
            \psit^{\Xbf^{k,\xbf}}_{k+1}
            ~-~
            \psit^{\Xbf^{k,\xbf, 0}}_{k+1} \1_{\{\Nt^k > 0 \}}
        \Big)
        ~\beta^{- \Nt^k}
        ~\prod_{j=1}^{\Nt^k}  \Wct^k_j.
    \ee
    We finally obtain the numerical algorithm of the path-dependent case:
     \be \label{eq:def_psib_p}
        \psit &:=& \psit^{x_0}_1.
    \ee

\subsubsection{The integrability and representation result}

    We notice that the algorithm in the path-dependent case is nothing else than an iterative algorithm of the Markovian case,
    as suggested by the PDEs \eqref{eq:PDE_n} in Remark \eqref{rem:PDE_n}.
    When the random variable $\psit$ in \eqref{eq:def_psib_p} is integrable, it is not surprising to obtain the representation $ \Vt_0 = \E \big[ \psit \big]$ as a consequence of Theorem \ref{theo:renorm_constant_vol}.
    However, because of the renormalization term 
    (i.e. $\big( \psit^{\Xbf^{k,\xbf}}_{k+1} - \psit^{\Xbf^{k,\xbf, 0}}_{k+1} \1_{\{\Nt^k > 0\}} \big)$ in \eqref{eq:def_psib}),
    the variance analysis becomes less obvious.
    We provide here a sufficient condition to ensure that $\psit$ has finite variance.

    \begin{Theorem} \label{theo:path_depend}
        Suppose that $\mu: [0,T] \x \R^{d \x n} \to \R^d$
        and $g: \R^{d \x n} \to \R$ are differentiable up to the order $n$,
        with bounded derivatives.
        Then
        \b*
            \E\big[ \big( \psit \big)^2 \big] ~<~ \infty
            &\mbox{and}&
            \Vt_0 ~:=~ \E \big[\psit \big].
        \e*
    \end{Theorem}

	We will prove the integrability result here, and leave the proof of the representation result $\Vt_0 := \E \big[\psit \big]$ in Section \ref{sec:Proofs}.
	As preparation, let us first provide two technical lemmas.
    Let $\pi = (0 = s_0 < s_1 < \cdots < s_m = T)$ be an arbitrary partition of the interval $[0,T]$, 
    $\mub:[0,T] \x \R^d \to \R^d$ a $\R^d-$valued function.
    We define $X^{\pi,x}$ by $X^{\pi,x}_0 := x$ and
    \be \label{eq:X_pi_x}
        X^{\pi,x}_{k+1}
        &:=&
        X^{\pi,x}_k
        ~+~
        \mub \big(s_k, X^{\pi,x}_k\big) \Delta s_{k+1}
        ~+~
        W_{s_{k+1}} - W_{s_k}.
    \ee
    Further, let $\varphi : \R^d \to \R$ be a smooth function, $\ell > 0$ and $i = (i_1, \cdots, i_{\ell}) \in \{1, \cdots, d\}^{\ell}$,
    we denote $\partial^{\ell}_{x,i} \varphi(x) := \partial^{\ell}_{x_{i_1} \cdots x_{i_{\ell}}} \varphi(x)$.

    \begin{Lemma} \label{lemm:X_diff}
        Suppose that $x \mapsto \mub(t,x)$ is differentiable up to order $n$ with uniformly bounded derivatives,
        and $X^{\pi,x}$ is defined by \eqref{eq:X_pi_x} with initial condition $X^{\pi,x}_0 = x$.
        Then $x \mapsto X^{\pi,x}_k$ is differentiable up to order $n$ and there is a constant $C$ independent of the partition $\pi$ such that
        \b*
            \max_{1 \le \ell \le n}
            ~\max_{i \in \{1, \cdots, d\}^{\ell}}
            ~\max_{0 \le k \le m}
            ~\big| \partial^{\ell}_{x,i} X^{\pi,x}_k \big|
            &\le&
            C.
        \e*
    \end{Lemma}
    \proof For simplicity, we consider the one dimensional $d=1$ case,
    while the multi-dimensional can be deduced by almost the same arguments.
    First, let $\ell = 1$, we have
    $$
        \partial_x X^{\pi,x}_{k+1}
        ~~=~ ~
        \partial_x X^{\pi,x}_k
        ~+~
        \partial_x \mub\big(s_k, X^{\pi,x}_k\big)  \partial_x X^{\pi,x}_k \Delta s_{k+1},
    $$
    which implies that
    $$
        \partial_x X^{\pi,x}_{k+1}
        ~~=~~
        \Pi_{j = 1}^{k+1} \Big( 1 + \partial_x \mub \big(s_k, X^{\pi,x}_k\big) \Delta s_{k+1} \Big).
    $$
    Since $\partial_x \mub(t,x)$ is uniformly bounded,
    it follows that $\partial_x X^{\pi, x}_k$ is bounded by some constant $C_1$ independent of $1 \le k \le m$ and the partition $\pi$.
    By induction, it is easy to deduce that for $\ell = 2, \cdots, n$,
    \b*
        \partial^{\ell}_{x^{\ell}} X^{\pi,x}_{k+1}
        ~~=~ ~
        \partial^{\ell}_{x^{\ell}} X^{\pi,x}_k
        ~+~
        P_{\ell} \big(\partial^i_{x^i} \mub(s_k, X_k^{\pi,x}), \partial^i_{x^i} X^{\pi,x}_k, i =1, \cdots, \ell -1 \big)
        \Delta s_{k+1},
    \e*
    where $P_{\ell}$ is a Polynomial on $\partial^i_{x^i} \mub(s_k, X_k^{\pi,x})$ and $\partial_{x^i}^i X^{\pi,x}_k$
    for $i =1, \cdots, \ell -1$,
    which is uniformly bounded by some constant independent of $k = 1, \cdots, m$ and the partition $\pi$.
    Hence $\partial^{\ell}_{x^{\ell}} X^{\pi,x}_k$ is also bounded by some constant $C_{\ell}$ independent of $k = 1, \cdots, m$ and the partition $\pi$.
    \qed

    \begin{Lemma} \label{lemm:estim_psi_k}
        Let $(\psit^{\xbf}_k)_{1 \le k \le n+1}$ be defined by \eqref{eq:def_psib}.
        Then for every $k= 2, \cdots, n+1$, and every $\xbf =(x_0, x_1, \cdots, x_{k-1}) \in \R^{d \x k}$,
        the map $x_{k-1} \mapsto \psit^{\xbf}_k$ has derivatives up to order $k-1$ and
        \be \label{eq:d_psib}
            \max_{1 \le \ell \le k-1}\Big| \partial^{\ell}_{x_{k-1},i}  \psit^{\xbf}_k \Big|
            &\le &
            C \prod_{j=k}^n (\Nt^j+1)^{j-1}.
        \ee
    \end{Lemma}
    \proof We will prove it by induction.
    First, let $k = n+1$, then $\psit^{\xbf}_{n+1} := g(x, x_1, \cdots, x_n)$ and hence
    $| \partial^{\ell}_{x_n} \psit^{\xbf}| \le C$ for some constant $C$ and for every $\ell = 1, \cdots, n$.

    Next, suppose that \eqref{eq:d_psib} holds true for $\psit^{\xbf}_{k+1}$,
    we know from \eqref{eq:def_psib} that
    \b*
        \psit^{\xbf}_k
        :=
        \Big(
            \psit^{\Xbf^{k,\xbf}}_{k+1}
            ~-~
            \psit^{\Xbf^{k,\xbf, 0}}_{k+1} \1_{\{\Nt^k > 0\}}
        \Big)
        \prod_{j=1}^{\Nt^k}  \frac{\mu_k(\Tt^k_j, \Xt^{k,\xbf}_j) - \mu_k(\Tt^k_{j-1}, \Xt^{k,\xbf}_{j-1})}{\beta \Delta \Tt^k_{j+1}} \cdot (\sigma_0^{\top})^{-1} \Delta \Wt^k_{j+1}.
    \e*
    Then using the estimation in Lemma \ref{lemm:X_diff}, we see that  \eqref{eq:d_psib} is also true for $\psit^{\xbf}_k$,
    and we hence conclude the proof.
    \qed

    \vspace{3mm}

    \noindent {\bf Proof of Theorem \ref{theo:path_depend} $\mathrm{(i)}$.}
    By Lemma \ref{lemm:estim_psi_k}, we know that
    $x \mapsto \psit_2^{x,x}$ is differentiable and in particular uniformly Lipschitz
    with coefficient bounded by $2 C \Pi_{j=2}^n (\Nt^j + 1)^{j-1}$.
    Then the definition of $ \psit_1^{x_0}$ falls into the Markovian case $n=1$, but with terminal condition
    $x \mapsto \psit^{x,x}_2$.
    Notice that $\Nt^k \le N_T$  has a Poisson distribution: $\P(N_T = m) = e^{-\beta T} \frac{(\beta T)^m}{m !}$.
    It follows that, for some constant $C > 0$,
    \b*
        ~\E \big[ \big| \psit_1^{x_0} \big|^2 \big]
        ~\le~
        \E \Big[
            C^{\Nt^k} 4 C^2 \prod_{j=2}^n (\Nt^j + 1)^{2 (j-1)}
        \Big]
        ~\le~
        \E \Big[
             4 C^2 C^{N_T} (N_T + 1)^{ n (n-1)}
        \Big]
        ~<~ \infty,
    \e*
    which implies that $\psit$ has finite variance.
	\qed

\section{Unbiased simulation of general SDEs}
\label{sec:general_vol}

	Let us now consider the SDE \eqref{eq:SDE} with general diffusion coefficient function, i.e. with drift and diffusion coefficients
	$\mu: [0,T] \x \R^d \to \R^d$ and $\sigma: [0,T] \x \R^d \to \M^d$:
	\b*
		X_0 = x_0,
		&\mbox{and}&
		dX_t
		~=~
		\mu \big(t, X_t \big) ~dt
		~+~
		\sigma \big(t, X_t\big) ~dW_t.
	\e*
	Our objective of study in this section is 
	$$
		V_0 ~=~ \E \big[ g(X_T) \big],
		~~\mbox{for some function}~ 
		g: \R^d \to \R.
	$$
	We will provide a representation result of $V_0$ in the same spirit of that in Section \ref{sec:constant_vol}.

    \begin{Remark}[Lamperti's transformation]
        \label{rem:Lamperti}

        We also notice that in some cases, the above SDE \eqref{eq:SDE} may be reduced
        to the constant diffusion coefficient case \eqref{eq:SDE3}, by the so-called the Lamperti transformation.

        \vspace{1mm}

        \noindent \rmi When $d = 1$ and $\sigma(t,x) > 0$,
        let us define a function $h : [0,T] \x \R \to \R$ by
        $$h(t, x) ~~:=~~ \int_0^x \frac{1}{\sigma(t, y)} dy.$$
        Notice that for fixed $t \in [0,T]$, 
        $x \mapsto h(t,x)$ is strictly increasing, we denote $h^{-1}(t, \cdot)$ its inverse function.
        Then by It\^o's formula, it is easy to obtain that $Y_t := h(t, X_t)$ satisfies the SDE
        \b*
            dY_t
            &=&
            \left( \partial_t h \big( t, h^{-1}(t,Y_t) \big)
                + \frac{\mu(t, h^{-1}(t,Y_t))}{\sigma(t, h^{-1}(t, Y_t))}
                - \frac{1}{2} \partial_x \sigma \big(t, h^{-1}(t, Y_t) \big)
            \right) dt
            ~+~
            dW_t    ,
        \e*
        whose diffusion coefficient is a constant as in SDE \eqref{eq:SDE3}.

        \vspace{1mm}

        \noindent \rmii When $d > 1$, $\sigma$ is non-degenerate and satisfies some further compatibility conditions,
        one can also obtain a similar transformation to reduce SDE \eqref{eq:SDE1}
        to the constant diffusion coefficient case.
    \end{Remark}

\subsection{An estimator of infinite variance for general SDEs }
\label{subsec:general_SDE}

	Let us impose the following conditions on coefficient functions $\mu$ and $\sigma$.
	\begin{Assumption} \label{assum:mu_sigma}
		The function $(\mu, \sigma)  : [0,T] \x \R^d \to \R^d \x \M^d$ and $a := \frac{1}{2} \sigma \sigma^{\top}: [0,T] \x \R^d \to \M^d$ are uniformly bounded,
		and are uniformly H\"older in the time variable, uniformly Lipschitz in the space variable, i.e. for some constant $L$,
		\be \label{eq:Lipschitz_coef2}
			\big| \big(\mu, \sigma, a \big) (t,x) - \big( \mu, \sigma, a\big) (s,y) \big|
			&\le&
			L \big( \sqrt{|t-s|} + \big| x-y \big| \big),
		\ee
		for all $ (t,x), (s,y) \in [0,T] \x \R^d$; 
		and $\sigma(t,x)$ is non-degenerate such that, for some constant $\eps_0 > 0$,
		\b*
			a(t,x) ~:=~ \frac{1}{2} \sigma \sigma^{\top}(t,x) ~\ge~ \eps_0 I_d,
			&&
			\forall (t,x) \in [0,T] \x \R^d.
		\e*
	\end{Assumption}
	
	Recall that $(T_k)_{k \ge 0}$ are defined by \eqref{eq:def_T_k} with a sequence of i.i.d. $\Ec(\beta)$-exponential random variables, and $W$ is a Brownian motion; the increment of the Brownian motion are defined by
	$\Delta W_{t_k} := W_{t_k} - W_{t_{k-1}}$, and $\Delta T_k := T_k - T_{k-1}$.
	As in \eqref{eq:REulerScheme}, we introduce $\Xh$ as solution of the Euler scheme on discrete grid 
	by $\Xh_0 := x_0$ and
	\be \label{eq:SDE_Euler_sol}
		\Xh_{T_{k+1}}
		~:=~
		\Xh_{T_k}
		+
		\mu\big(T_k, \Xh_{T_k} \big) \Delta T_{k+1}
		+
		\sigma \big(T_k, \Xh_{T_k}\big) \Delta W_{T_{k+1}}
		~~~k = 0, \cdots, N_T.
	\ee
	We then introduce a representation formula by
	\be \label{eq:def_psih2}
		\psih
		&:=&
		e^{\beta T} ~
		\big[
		g \big(\Xh_T \big)
		-
		g \big( \Xh_{T_{N_T}} \big) \1_{\{N_T > 0\}}
		\big]
		~\beta^{-N_T}
		~\prod_{k=1}^{N_T} \Big( \Wcb^1_k + \Wcb^2_k \Big),
	\ee
	where, for each $k = 1, \cdots, N_T$,
	\b*
		\Wcb^1_k
		&:=&
		\big[\mu(T_k, \Xh_{T_k}) - \mu(T_{k-1}, \Xh_{T_{k-1}})
		\big]
		\cdot~
		\frac{
		\big(\sigma^{\top} (T_{k}, \Xh_{T_k})\big)^{-1} \Delta W_{T_{k+1}}}
		{ \Delta T_{k+1}},
  	\e*
	and
	\be \label{eq:Wcb2}
		\Wcb^2_k
		&:=&
		\big[
			a \big(T_k, \Xh_{T_k} \big)
			- a \big(T_{k-1}, \Xh_{T_{k-1}} \big)
		\big] \nonumber \\
		&&~~
		:
		\Big[ \big(\sigma^{\top}(T_{k}, \Xh_{T_k}) \big)^{-1}
		~\frac{\Delta W_{T_{k+1}} \Delta W^{\top}_{T_{k+1}} - \Delta T_{k+1} I_d }{ \Delta T^2_{k+1}}
		~\sigma(T_{k}, \Xh_{T_k})^{-1}
		\Big].
	\ee

	\begin{Theorem} \label{theo:renorm_general_vol}
		Suppose that Assumption \ref{assum:mu_sigma} holds true, and $g$ is Lipschitz.
		Then
		\b*
			\E \big[ \big| \psih \big| \big]
			~<~
			\infty
			&\mbox{and}&
			V_0 = \E [ \psih].
		\e*
	\end{Theorem}
	\proof
	\rmi Consider the random vectors $\xi_k^1 := \frac{\Delta W_{T_k}}{\sqrt{\Delta T_k}}$ and $\xi_k^2:= \frac{\Delta W_{T_k} \Delta W^{\top}_{T_k} - \Delta T_k I_d}{\Delta T_k}$, for all $k = 1, \cdots, N_T +1$,
	which are independent of $\Delta T_k$ conditional on $\{\Delta T_k >0\}=\{N_T \ge k-1\}$, and which have finite second order moment.
	Notice that $\mu(t,x)$ and $a(t,x)$ are uniformly bounded, and $1/2-$H\"older-continuous in $t$ and Lipschitz in $x$, and $\sigma$ is uniformly bounded from below above zero. Then, for each $k=1, \cdots, N_T$,
	$$
		\big| \Wcb^1_k \big|
		\le
		C \big( \sqrt{\Delta T_k} + \big| \Xh_{T_k} - \Xh_{T_{k-1}} \big| \big) \Big| \frac{\Delta W_{T_{k+1}}}{\Delta T_{k+1}} \Big|
		\le
		C \big( 1 + \sqrt{\Delta T_k} + |\xi^1_k| \big) |\xi^1_{k+1}| \sqrt{\frac{\Delta T_k}{\Delta T_{k+1}}},
	$$
	where the constant $C > 0$ may vary from term by term but is uniformly bounded for all $k$.
	Similarly, one obtains that
	$$
		\big| \Wcb^2_k \big|
		~\le~
		C \big( 1 + \sqrt{\Delta T_k} + |\xi^1_k| \big) |\xi^2_{k+1}| \sqrt{\frac{\Delta T_k}{\Delta T_{k+1}}} \frac{1}{\sqrt{\Delta T_{k+1}}}.
	$$
	As $\Delta T_k \le T$, it follows that
	$$
		\big| \Wcb^1_k \big| + \big| \Wcb^2_k \big|
		~\le~
		C \big( 1 + |\xi_k^1| \big) \big( |\xi^1_{k+1}| + |\xi^2_{k+1}| \big) \sqrt{\frac{\Delta T_k}{\Delta T_{k+1}}} \frac{1}{\sqrt{\Delta T_{k+1}}},
	$$
	for some constant $C > 0$ independent of $k$. In addition, we have by the Lipschitz condition on $g$ that
	$$ 
		\E \Big[ \big | g(\Xh_T) - g(\Xh_{T_{N_T}}) \big|  ~\Big| \Delta T_{N_T+1} \Big] 
		~<~
		C \sqrt{\Delta T_{N_T+1}}.
	$$
	Then, it follows from the expression of $\psih$ in \eqref{eq:def_psih2} that
	\b*
		\E\big[ \big| \psih \big| \big]
		&\le&
		C \E \Big[
		\prod_{k=1}^{N_T} \frac{C}{\sqrt{\Delta T_{k+1}}}
		\1_{\{ N_T \ge 1\}} \Big]
		+
		C \E \Big[  \big| g(\Xh_T) \big|~\1_{\{ N_T = 0\}} \Big] \\
		&\le&
		C \E \Big[
		\prod_{k=1}^{N_T} \frac{C}{\sqrt{\Delta T_{k+1}}} \Big]
		+
		C \E \Big[  \big| g(x_0 + \mu(0,x_0) T + \sigma(0,x_0) W_T ) \big| \Big] \\
	\e*
	for some constant $C > 0$, where we have also used the independence of the $\xi_k^i$'s and their the boundedness of their second order moments.
	The integrability of $\psih$ is now a direct consequence of Lemma \ref{lemm:OrderStat_p}.

	\vspace{1mm}

	\noindent \rmii The proof of the equality $V_0 = \E [ \psih]$ will be completed in Section \ref{sec:Proofs}.
	\qed

	\vspace{2mm}
	
	To conclude, we notice that the variable $\psih$ is of order $\Pi_{k=1}^{N_T} 1/ \sqrt{\Delta T_{k+1}}$ in general cases,
        and the latter is integrable but of infinite variance.
        Therefore, $\psih$ is not a good estimator for Monte-Carlo method.
	Nevertheless, it should still have some theoretical value as an alternative representation formula obtained by 
	Bally and Kohatsu-Higa \cite[Section 6.1]{Bally}.

\subsection{An estimator for one-dimensional driftless SDE}

	To overcome the problem of variance explosion of the estimator \eqref{eq:def_psih2}, 
	we will consider the higher order approximation $\Xh$ of $X$,
	and obtain an estimator of finite variance for the one dimensional ($d = 1$) driftless SDE of form
	\be
		X_0 = x_0, ~~~~dX_t ~=~ \sigma(t, X_t) ~dW_t,
	\ee
	Our objective is  to compute
	\b*
		V_0 &:=& \E \big[ g(X_T) \big], ~~~\mbox{for some function}~g:\R \to \R.
	\e*

	Recall $(T_k)_{k \ge 0}$ has been introduced in \eqref{eq:def_T_k} from a sequence of i.i.d. exponential random variables, independent of the Brownian motion $W$.
	We next define $\Xh$ by $\Xh_0 = x_0$,
	\be \label{eq:dXh_2}
		d \Xh_t
		~=~
		\Big( \sigma(T_k, \Xh_{T_k}) + \partial_x \sigma(T_k, \Xh_{T_k}) \big( \Xh_t  - \Xh_{T_k}\big) \Big) dW_t,
		&\mbox{on}~
		[T_k, T_{k+1}],
	\ee
    for $k = 0, 1, \cdots, N_T$.
    By denoting 
    \be \label{eq:def_c1c2k}
    	c_1^k := \sigma(T_k, \Xh_{T_k}) -  \partial_x \sigma(T_k, \Xh_{T_k}) \Xh_{T_k}
    	&\mbox{and}&
	c_2^k := \partial_x \sigma(T_k, \Xh_{T_k}),
    \ee
    then the above linear SDE \eqref{eq:dXh_2} has an explicit solution which is given by
    \be \label{eq:linSDE_sol1}
        \Xh_{T_{k+1}} = \Xh_{T_k} + \sigma(T_k, \Xh_{T_k})  \Delta W_{T_{k+1}},
        ~~~~~\mbox{if}~~
        c_2^k = 0,
    \ee
    and
    \be \label{eq:linSDE_sol2}
        ~\Xh_{T_{k+1}}
        &=&
        -~ \frac{c^k_1}{c^k_2}
        ~+~ \frac{c^k_1}{c^k_2} \exp \Big( - \frac{(c^k_2)^2}{2} \Delta T_{k+1} + c^k_2 \Delta W_{T_{k+1}} \Big) \nonumber \\
        &&
        +~ \Xh_{T_k} \exp \Big(- \frac{(c^k_2)^2}{2} \Delta T_{k+1} + c^k_2 \Delta W_{T_{k+1}} \Big),
        ~~~~~~~~~\mbox{if}~ c_2^k \neq 0.
    \ee
	We then define $\psih$ by
    \be \label{eq:psih2}
        \psih
        &:=&
         e^{\beta T}
        \Big[
            g(\Xh_T) -  g(\Xh_{T_{N_T}})  \1_{\{N_T > 0\}}
        \Big]
        ~\beta^{-N_T}
        ~\prod_{k=1}^{N_T}
        \Wcb^2_k,
    \ee
    where the automatic differentiation  weight is given by (see Lemma \ref{lemm:MalliavinW_LogNorm} below)
    \be \label{eq:MWeight2}
        \Wcb^2_k
        ~:=~
                \frac{a(T_k, \Xh_{T_k}) - \tilde a_k}{ 2 a(T_k, \Xh_{T_k})}
            \Big(-~ \partial_x \sigma(T_k, \Xh_{T_k}) \frac{\Delta W_{T_{k+1}}}{\Delta T_{k+1}}
            ~+~ \frac{\Delta W_{T_{k+1}}^2 -\Delta T_{k+1}}{\Delta T_{k+1}^2} \Big),
    \ee
    with $a(\cdot) := \frac{1}{2} \sigma^2(\cdot)$, $\tilde a_k := \frac{1}{2} \tilde \sigma_k^2$ and
    $\tilde \sigma_k  := \sigma(T_{k-1}, \Xh_{T_{k-1}})
    + \partial_x \sigma(T_{k-1}, \Xh_{T_{k-1}}) (\Xh_{T_k} - \Xh_{T_{k-1}})$.

	Similarly to the discussion at the end of Section \ref{subsec:general_SDE} 
	(see also Remark \ref{rem:generalTHM} below), 
	the variable $\psih$ in \eqref{eq:psih2} is integrable but of infinite variance in general.
	To make the variance finite, we introduce an alternative estimator using an antithetic variable.
	Let $\Xh^{-}_T$ be an antithetic variable of $\Xh_T$ defined by
	\b*
		\Xh^-_T ~:=~
		\Xh_{T_{N_T}}
		- \sigma(T_{N_T}, \Xh_{T_{N_T}}) \Delta W_{T_{N_T}},
		~~~~\mbox{if}~ c_2^{N_T} = 0,
	\e*
    and
    \b*
        ~\Xh^-_T
        &=&
        -~ \frac{c^{N_T}_1}{c^{N_T}_2}
        ~+~ \frac{c^{N_T}_1}{c^{N_T}_2} \exp \Big( - \frac{(c^{N_T}_2)^2}{2} \Delta T_{N_T+1}
                - c^{N_T}_2 \Delta W_{T_{N_T+1}} \Big) \\
        &&
        +~ \Xh_{T_{N_T}} \exp \Big(- \frac{(c^{N_T}_2)^2}{2} \Delta T_{{N_T}+1}
                - c^{N_T}_2 \Delta W_{T_{N_T+1}} \Big),
        ~~~~~~~~~\mbox{if}~ c_2^{N_T} \neq 0.
    \e*
    Denote $\Wcb_k^- := \Wcb^2_k$ for $k = 1, \cdots, N_T -1$ and
    \b*
        \Wcb^-_{N_T} 
        &:=&
        \frac{a(T_{N_T}, \Xh_{N_T}) - \tilde a_{N_T}}{ 2 a(T_{N_T}, \Xh_{N_T})}
	\Big(
            \partial_x \sigma(T_{N_T}, \Xh_{N_T}) \frac{\Delta W_{T_{N_T+1}}}{\Delta T_{{N_T}+1}}
            +
            \frac{\Delta W_{T_{N_T+1}}^2 -\Delta T_{{N_T}+1}}{\Delta T_{{N_T}+1}^2}
        \Big).
    \e*
    We then introduce
    \be \label{eq:psib2}
        \psib := \frac{\psih + \psih^-}{2}
        ~~\mbox{with}~~
        \psih^-
        :=
         e^{\beta T}
        \Big[
            g(\Xh^-_T) -  g(\Xh_{T_{N_T}})  \1_{\{N_T > 0\}}
        \Big]
        ~\beta^{-N_T}
        ~\prod_{k=1}^{N_T}
        \Wcb_k^-.
    \ee
    Notice that the Brownian motion is symmetric, thus $\psih^-$ has exactly the same distribution as $\psih$,
    and it serves as an antithetic variable.

	\begin{Assumption} \label{assum:1d_driftless}
		The diffusion coefficient  $\sigma(\cdot)$ satisfies 
		$\sigma(t,x) \ge \eps > 0$ for all $(t,x) \in [0,T] \x \R$,
		$\sigma(t,x)$ is bounded and Lipschitz in $(t,x)$,
		$\partial_x \sigma(t,x)$ is bounded continuous in $(t,x)$ and uniformly Lipschitz in $x$.
		Further, the terminal condition function $g(\cdot) \in C_b^2(\R)$.
	\end{Assumption}

    \begin{Theorem} \label{theo:ExactSimulation2}
        Suppose that Assumption \ref{assum:1d_driftless} holds true.
        Then
        \be
            \E \big[ \big| \psih \big| \big] 
            ~ +~
            \E \big[ \big| \psib \big|^2 \big] ~<~\infty;
            &\mbox{and}&
            V_0 ~=~ \E\big[ \psih \big] ~=~ \E \big[ \psib \big].
        \ee
    \end{Theorem}

	We will complete the proof in Section \ref{subsec:proof_0drift1d}.

    \begin{Remark}
	\rmi As $\psib$ has  finite variance, we may use the representation of Theorem \ref{theo:ExactSimulation2}
	to built an unbiased Monte-Carlo estimator of $V_0$.
	However, given the assumed regularity conditions, and the restriction to the one-dimensional setting,
	such a Monte-Carlo approximation is not competitive with the corresponding PDE based approximation methods.
	However, we believe that the present methodology is open to potential improvements, and we hope to improve our results in some future work so as to address the higher dimensions.
	
	\noindent \rmii For a general SDE with drift function and/or $d \ge 1$, we can also consider a similar choice of $(\muh, \sigmah)$,
        which leads to $\muh(t,x) = c_1 + c_2 x$ and $\sigmah(t,x) = c_3 + c_4 x$ and a linear SDE
        \be \label{eq:Ge_Lin_SDE}
            d\Xh_t
            &=&
            \big(c_1 + c_2 \Xh_t \big) dt
            ~+~
            \big(c_3 + c_4 \Xh_t \big) dW_t,
        \ee
        where $c_1 \in \R^d$, $c_2, c_3 \in \M^d$ and $c_4$ is linear operator from $\R^d$ to $\M^d$.
        However, to the best of our knowledge,
        the exact simulation of linear SDE \eqref{eq:Ge_Lin_SDE} in high dimensional case,
        as well as the associated automatic differentiation (Malliavin) weight as in \eqref{eq:MWeight2} (see also Lemma \ref{lemm:MalliavinW_LogNorm} below),
        is still an open question.
    \end{Remark}

\section{Numerical examples}
\label{sec:Numerics}

	Notice that our estimator $\psih$ given by \eqref{eq:def_psih1} (resp. $\psit$ given by \eqref{eq:def_psib} and \eqref{eq:def_psib_p}) is an unbiased estimator for $V_0$ in \eqref{eq:def_V0} (resp. $\Vt_0$ in \eqref{eq:def_Vt0}).
	Then the error analysis of the Monte-Carlo approximation reduces to the statistical error.
	Hence the computation cost to achieve the accuracy $O(\eps)$ for the approximation of $V_0$ (resp. $\Vt_0$) is of order $O(\eps^{-2})$,
	thus avoiding of the dependence on the discretization error.
	
	By combining different level of simulations, the MultiLevel Monte Carlo (MLMC) method proposed by Giles \cite{Giles} achieves a computation cost of order $O(\eps^{-2} (\log \eps)^2)$ or $O(\eps^{-2})$ depending on the strong discretization error rate.
	In particular, by considering a randomization of the level, Rhee and Glynn \cite{RheeGlynn} obtained an unbiased estimator.
	In the following, we provide some numerical results and comparisons between our unbiased simulation method with the Euler based MLMC method proposed by \cite{Giles}.

\subsection{Two one-dimensional SDEs}

	Let $W$ be a one-dimensional standard Brownian motion, we consider the SDE given by
	$$
		S_0 = 1,
		~~~~ dS_t ~=~ 0.1 \big(\sqrt{M \wedge S_t} - 1 \big) S_t dt ~+~ \frac{1}{2} S_t dW_t,
	$$
	where $M$ is a large constant introduced in order to guarantee the Lipschitz property of the drift coefficient (in our numerical implementation, we have observed that the value of $M$ is not relevant for large $M$, and that the numerical finding are not changed by taking $M=\infty$; this hints that our results may be extended beyond the case of Lipschitz coefficients). Applying Lemperti's transformation $X_t := \log(S_t)$, we reduce the above SDE to the constant diffusion coefficient case, in form of \eqref{eq:SDE3},
	\be \label{eq:SDE_num11}
		X_0 = 0,
		~~~ dX_t ~=~ \Big( 0.1 \big( \sqrt{M \wedge e^{X_t}} - 1 \big) - 1/8 \Big) dt ~+~ \frac{1}{2} dW_t.
	\ee
	We implement our unbiased simulation method for the two following expectations:
	\be \label{eq:def_V0_num1}
		V_0 ~:=~ \E \big[ (S_T - K)_+ \big],
		~~~\mbox{and}~~
		\Vt_0 ~:=~ \E \Big[ \Big( \frac{1}{n} \sum_{k=1}^n S_{t_k} - K \Big)_+ \Big],
	\ee
	where we choose $K = 1$, $T=1$, $n = 10$ and $t_k := \frac{k}{n} T$.
	Notice that the path-dependent example does not satisfy the differentiability sufficient condition in Theorem \ref{theo:path_depend}.
	However, our numerical findings do not show any numerical difficulty in the present setting.

	Using different numbers $N$ of simulations, we obtain the standard deviation as (statistical) error of our estimator.
	Next, using the errors obtained by our unbiased simulation method,
	we implement the MLMC algorithm in Section 5 of Giles \cite{Giles}, and we compare the computation time (in second) of the two methods.
	More precisely, the statistical error of the unbiased simulation method is given by $\sqrt{\mathrm{Var}[\psih] / N}$, where $\mathrm{Var}[\psih] $ denotes the estimated variance of $\psih$. 
	For the implementation of  MLMC, we choose $M=4$, $N_L = 10^4$ and use equation (10) in \cite{Giles} as criteria to stop the loop in MLMC (see more details in Section 5 of \cite{Giles} for the meaning of $M$ and $N_L$).

\begin{table}
\begin{center}
\begin{tabular}{|c||c|c|c|}
\hline
~   &  Mean value  & Statistical error & Computation time \\ \hline
US ($N = 10^5$) & 0.204864 & 0.00140709 & 0.016814 \\ \hline
MLMC & 0.204993 & 0.000949166 & 0.032017 \\ \hline \hline
US ($N = 10^6$) & 0.205396 & 0.000444462 & 0.171835 \\ \hline
MLMC & 0.205602 & 0.000308634 & 0.234526 \\ \hline \hline
US ($N = 10^7$) & 0.20552 & 0.000142554 & 1.63013 \\ \hline
MLMC & 0.205648 & 0.0001 & 1.96197 \\ \hline \hline
US ($N = 10^8$) & 0.205641 & 4.52282e-05 & 16.2189 \\ \hline
MLMC & 0.205638 & 3.18855e-05 & 18.3833 \\ \hline
\end{tabular}
\end{center}
\caption{Numerical results for $V_0$ in \eqref{eq:def_V0_num1} (case $d=1$), US denotes our unbiased simulation algorithm with $\beta = 0.1$,
the computation times are expressed in second. }
\label{Tab:SDE1_1}
\end{table}

\begin{table}
\begin{center}
\begin{tabular}{|c||c|c|c|}
\hline
~   &  Mean value  & Statistical error & Computation time \\ \hline
US ($N = 10^5$) & 0.127032 & 0.000762635 & 0.144998 \\ \hline
MLMC & 0.127053 & 0.000536248 & 0.323337 \\ \hline \hline
US ($N = 10^6$) & 0.126363 & 0.000241231 & 1.40843 \\ \hline
MLMC & 0.126747 & 0.000169842 & 1.8194 \\ \hline \hline
US ($N = 10^7$) & 0.126703 & 7.6418e-05 & 13.9005 \\ \hline
MLMC & 0.126643 & 5.37691e-05 & 16.7499 \\ \hline
\end{tabular}
\end{center}
\caption{Numerical results for $\Vt_0$ in \eqref{eq:def_V0_num1} (case $d=1$), US denotes our unbiased simulation algorithm with $\beta = 0.05$,
the computation times are expressed in second. }
\label{Tab:SDE1_2}
\end{table}

\begin{figure}\begin{center}
\includegraphics[width=14cm,height=9cm]{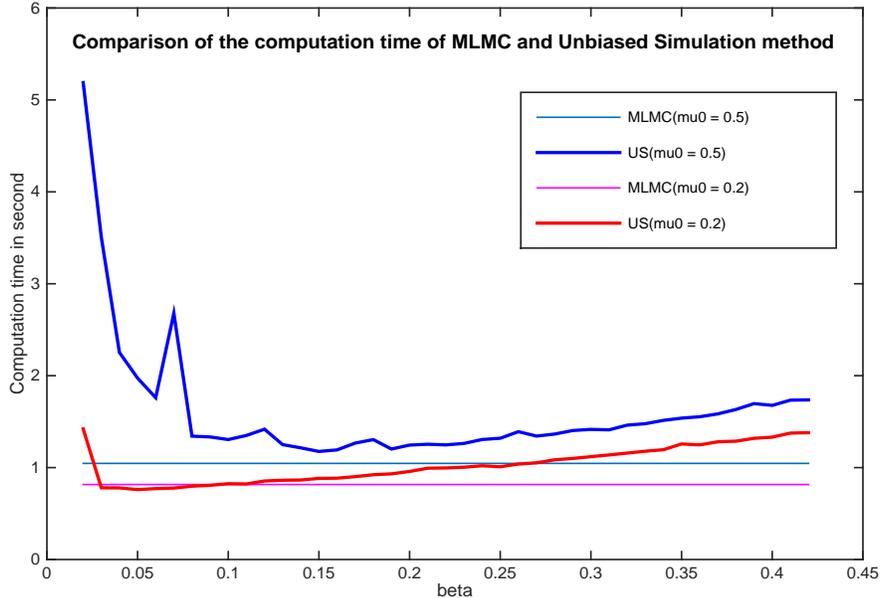}
\end{center}
\caption{Comparison of the computation time of MLMC method and unbiased simulation method for problem \eqref{eq:def_V0_num11}, with the same given error.} \label{Fig0}
\end{figure}

	The numerical results are given in Tables \ref{Tab:SDE1_1} and \ref{Tab:SDE1_2}.
	We observe that with the same Monte-Carlo error, both methods have very close performance.
	In the present particular example, the computational time of our methods is slightly smaller.
	However, the conclusion may change depending on the nature of the example.
	Let us consider the problem
	\be \label{eq:def_V0_num11}
		V_0 ~:=~ \E\big[ \sin(X_T) \big],
	\ee
	where $X$ is defined by SDE, for some constant $\mu_0 \in \R$,
	\b*
		X_0 = 0, && dX_t ~=~ \mu_0 \cos(X_t) dt ~+~ \frac{1}{2} dW_t.
	\e*
	We implement the MLMC algorithm and our unbiased simulation method with different value of $\beta$,
	but with a given fixed error $\varepsilon = 0.0002$.
	The two methods provide very close estimation of value $V_0$, so we give a comparison on the computation time in Figure \ref{Fig0}.
	We can observe that $\beta$ in the unbiased simulation method should not be too big nor too small, to minimize the computation effort.
	When $\mu_0 = 0.2$, the computation time of MLMC method is slightly longer than the US method with $\beta \approx 0.05$.
	However, when $\mu_0 = 0.5$, the computation time MLMC method is always smaller than the US method for any choice of $\beta > 0$. This shows that, in the context of the present example, the performance of our unbiased simulation method is of the order of that of the multilevel Monte Carlo method.

\subsection{A multi-dimensional SDE}

	We next consider a d-dimensional SDE with $d=4$.
	Let $W = (W^1, \cdots, W^4)^{\top}$ be a $4$-dimensional standard Brownian motion, and 
	$\sigma_0$ the $4 \x 4$ be the lower triangular matrix such that
	\[
		\sigma_0 \sigma_0^{\top}
		~=~
		\left( \begin{array}{cccc}
		1 & 1/2 & 1/2 & 1/2 \\
		1/2 & 1 & 1/2 & 1/2 \\
		1/2 & 1/2 & 1 & 1/2 \\
		1/2 & 1/2 & 1/2 & 1
		 \end{array} \right).
	\] 
	We consider the SDE
	\b*
		dX_t ~=~ \mu(t, X_t)  dt ~+~   \sigma_0 d W_t,
		~~~X^i_0 = 0, ~i = 1,\cdots, 4,
	\e*
	with drift function $\mu(t,x) = (\mu_i(t,x), ~i =1, \cdots, 4)$ be given by 
	$\mu_i(t,x_1, \cdots, x_4) = 0.1 \Big( \sqrt{ \frac{3}{4} \exp(x_i) + \frac{1}{4} \overline{\exp(x)} } - 1 \Big) - \frac{1}{8}$,
	where $\overline{\exp(x)} := (e^{x_1} + \cdots + e^{x_4})/4$.
	We then consider two problems:
	\be \label{eq:def_V0_num2}
		V_0 := \E \Big[ \Big(\frac{1}{4} \sum_{i=1}^4 e^{M \wedge X^i_T} - K \Big)_+ \Big],
		~~\mbox{and}~
		\Vt_0 :=\E \Big[ \Big( \frac{1}{4n} \sum_{k=1}^n \sum_{i=1}^4 e^{M \wedge X^i_{t_k}} - K \Big)_+ \Big],
	\ee
	where we choose $K = 1$, $T=1$, $n = 10$ and $t_k := \frac{k}{n} T$ and $M$ is a large number so as to ensure that the terminal condition is Lipschitz.
	As in the one-dimensional case, we implement our unbiased simulation method using different sample sizes $N$.
	Then, we use the errors, obtained from our unbiased simulation method, in the MLMC algorithm in Section 5 of Giles \cite{Giles}, and we compare the computation time (in second) of the two methods.

	The numerical results are given in Tables \ref{Tab:SDE2_1} and \ref{Tab:SDE2_2}.
	We observe that both methods have very similar performance, with a slightly small advantage for our method.
	However, similar to the one-dimensional case, the MLMC algorithm could be better in other examples.

\begin{table}
\begin{center}
\begin{tabular}{|c||c|c|c|}
\hline
~   &  Mean value  & Statistical error & Computation time \\ \hline
US ($N = 10^5$) & 0.739374 & 0.00921078 & 0.109151 \\ \hline
MLMC & 0.732707 & 0.00568921 & 0.136884 \\ \hline \hline
US ($N = 10^6$) & 0.735745 & 0.00239613 & 1.06639 \\ \hline
MLMC & 0.733539 & 0.00176862 & 1.15886 \\ \hline \hline
US ($N = 10^7$) & 0.73659 & 0.000831597 & 10.6957 \\ \hline
MLMC & 0.737087 & 0.000578058 & 12.171 \\ \hline 
\end{tabular}
\end{center}
\caption{Numerical results for $V_0$ in \eqref{eq:def_V0_num2} (case $d=4$), US denotes our unbiased simulation algorithm with $\beta = 0.5$,
the computation times are expressed in second. }
\label{Tab:SDE2_1}
\end{table}

\begin{table}
\begin{center}
\begin{tabular}{|c||c|c|c|}
\hline
~   &  Mean value  & Statistical error & Computation time \\ \hline
US ($N = 10^5$) & 0.382186 & 0.00247547 & 0.769847 \\ \hline
MLMC & 0.381071 & 0.00167112 & 2.07589 \\ \hline \hline
US ($N = 10^6$) & 0.382846 & 0.000762393 & 7.65796 \\ \hline
MLMC & 0.383107 & 0.000535905 & 10.8444 \\ \hline \hline
US ($N = 10^7$) & 0.383282 & 0.000244861 & 85.0265 \\ \hline
MLMC & 0.383653 & 0.00017245 & 104.223 \\ \hline
\end{tabular}
\end{center}
\caption{Numerical results for $\Vt_0$ in \eqref{eq:def_V0_num2} (case $d=4$), US denotes our unbiased simulation algorithm with $\beta = 0.05$,
the computation times are expressed in second. }
\label{Tab:SDE2_2}
\end{table}

\subsection{A one-dimensional driftless SDE}

	Finally, we provide an example of a one-dimensional driftless SDE.
	We recall that under the assumed regularity in Theorem \ref{theo:ExactSimulation2}, we are not expecting our method to be competitive with the PDE based approximations.
	Instead, our objective  is to study numerically the performance of the estimator \eqref{eq:psib2}.

	Let us consider the SDE
	\be \label{eq:SDE_num3}
		X_0 = 1, ~~~ dX_t = \frac{2 \sigma}{	1 + X^2_t} dW_t,
	\ee
	and we aim to compute
	\be \label{eq:def_V0_num3}
		V_0 ~=~ \E \big[ \big(X_T - K \big)_+ \big].
	\ee
	We implemented the simple Euler scheme with time step $\Delta t = 1/10$ and simulation number $N = 10^6$,
	and next the unbiased simulation method \eqref{eq:psib2} with $\beta = 0.1$ and simulation number $N = 10^6$,
	and then the MLMC scheme using the statistical error obtained from the unbiased simulation method.
	The results are given in Table \ref{Tab:num3}, and we can observe that all three methods provides very similar estimation of $V_0$.
	In particular, the unbiased simulation method has a significant advantage.

\begin{table}
\begin{center}
\begin{tabular}{|c||c|c|c|}
\hline
~   &  Mean value  & Statistical error & Computation time \\ \hline
Euler scheme & 0.161483 & 0.000196733 & 0.570541 \\ \hline
US &0.160362 &9.34729e-05  & 0.201904 \\ \hline 
MLMC & 0.16057 & 6.61696e-05 & 6.65799 \\ \hline
\end{tabular}
\end{center}
\caption{Numerical results for $V_0$ in \eqref{eq:def_V0_num3} (case $d=1$), 
US denotes the unbiased simulation algorithm \eqref{eq:psib2} with $\beta = 0.1$,
the computation times are expressed in second. 
	Notice also that the unbiased algorithm \eqref{eq:psib2} contains implicitly an antithetic variance reduction,
	which makes its statistical even error smaller than that of the Eurler scheme.
}
\label{Tab:num3}
\end{table}

\section{Proofs}
\label{sec:Proofs}

\subsection{A toy example}
\label{subsec:toy_proof}

	Before completing the technical part of the proofs for Theorems \ref{theo:renorm_constant_vol}, \ref{theo:path_depend} and \ref{theo:renorm_general_vol}, 
	we would like to illustrate the main idea by studying a simplified example in the one dimensional case 
	with unit diffusion:
	$$
		X_0 = x_0,
		~~~~~
		dX_t ~=~ \mu(t, X_t) dt ~+~ dW_t.
	$$
	
	Let $b \in \R$ and $\beta > 0$, we define a sequence of i.i.d. random variable $(\tau_k)_{k\ge 1}$ of distribution $\Ec(\beta)$.
	Then let $(T_k)_{k \ge 1}$, $\Delta T_{k+1} := T_{k+1} - T_k$ and $\Delta W_{k+1} := W_{T_{k+1}} - W_{T_k}$ be defined in and below \eqref{eq:def_T_k},
	we introduce $\Xh$ by 
	\b*
		\Xh_t &:=& x_0 ~+~ b t ~+~  W_t.
	\e*
	and then define
	\be \label{eq:psit_vol_constant}
		\psi
		&=&
		e^{\beta T} ~g \big(\Xh_T \big) 
			\prod_{k=1}^{N_T} 
			~\frac{ \big( \mu(T_k, \Xh_{T_k}) - b \big)  \Delta W_{T_{k+1}}}
			{\beta \Delta T_{k+1}}.
	\ee
		
	\begin{Proposition} \label{prop:brut_algo}
		Let $\mu(\cdot, \cdot)$ and $g(\cdot)$ be both bounded smooth functions in $C_b^2$.
		Then for all constants $\beta > 0$ and $\beta > 0$,
		one has
		$$
			\E \big[ \big| \psi \big| \big] ~<~ \infty
			~~~\mbox{and}~~~
			\E \big[ g(X_T) \big] ~=~ \E \big[ \psi \big].
		$$
	\end{Proposition}
	\proof Since $\mu$ and $g$ are uniformly bounded,
	and for some constant $C > 0$,
	the conditional expectation $\E_{\Delta T_{k+1}} [ | \Delta W_{T_{k+1}}|] \le C/\sqrt{\Delta T_{k+1}} $,
	then by Lemma \ref{lemm:OrderStat_p}, it is obvious that $\psi$ is integral. 
	Then it is enough to prove that $\E \big[ g(X_T) \big] = \E \big[ \psi \big]$.
	In preparation, let us introduce
	\b*
		\psi_n
		&=&
		e^{\beta T_{n+1}}
		\prod_{k=1}^{N_T \wedge n} \frac{(\mu(T_k, \Xh_{T_k}) - b) \Delta W_{T_{k+1}}}{\beta \Delta T_{k+1}}\\
		&&~~\Big( g \big(\Xh_T \big) \1_{N_T \le n} 
		+
		\big(\frac{\mu - b}{\beta}  \partial_x u \big) \big(T_{n+1},\Xh_{T_{n+1}}\big) \1_{N_T > n}
		\Big).
	\e*
	for all $n \ge 0$, with the convention $\prod_{k=1}^0 \equiv 1$.
	It is clear that $(\psi_n)_{n \ge 0}$ are all integrable by  Lemma \ref{lemm:OrderStat_p}.

	\vspace{1mm}
	
	\noindent \rmi Notice that $\mu$ and $g$ are both smooth functions, then by Feynmann-Kac formula,
	we know $\E[g(X_T)] = u(0,x_0)$,
	where $u \in C_b^{\infty}([0,T] \x \R)$ is a smooth function of PDE
	$$
		\partial_t u (t,x) ~+~ \frac{1}{2} \partial^2_{xx} u (t,x)  ~+~ \mu(t,x) \partial_x u (t,x) ~=~ 0,
		~~~\mbox{for all}~~(t,x) \in [0,T) \x \R,
	$$
	with terminal condition $u(T,x) = g(x)$.
	Rewriting the above PDE in the following equivalent way:
	$$
		\partial_t u(t,x) ~+~ b \partial_x u(t,x) ~+~  \frac{1}{2} \partial^2_{xx} u(t,x)
		~+~ \big(\mu(t,x) - b \big) \partial_x u(t,x) = 0,
	$$
	it follows from the  Feynmann-Kac formula that
	\be \label{eq:Expr1_toy}
		u(0,x_0) 
		\!\!&=&\!\!
		\E \Big[ g(\Xh_T) + \int_0^T \big(\mu\big(t, \Xh_t\big) - b \big) \partial_x u \big(t,\Xh_t\big) dt \Big] \nonumber \\
		\!\!&=&\!\!
		\E \Big[ e^{\beta T} g \big(\Xh_T \big)\1_{\{T_1 \ge T\}} + \frac{e^{\beta T_1}}{\beta} \big(\mu\big(T_1, \Xh_{T_1}\big) - b \big)  \partial_x u \big(T_1,\Xh_{T_1}\big) \1_{\{T_1 < T\}} 
		\Big]~~ ~~\\
		\!\!&=&\!\! \E\big[ \psi_0 \big], \nonumber
	\ee
	where the second equality follows from the fact that $T_1 = T \wedge \tau_1$, and $\tau_1$ is a random variable independent of $\Xh$, 
	with density function $\beta e^{- \beta t} \1_{\{t \ge 0\}}$.

	\vspace{1mm}

	\noindent \rmii Next, notice that for any $t$ and bounded continuous function $\phi_0$, one has by integration by parts that
	\be \label{eq:DeriveGaussian_toy}
		\partial_x \E\big[ \phi_0(x + b t + W_t) \big] ~=~ \E \Big[ \phi_0(x + b t + W_t) \frac{W_t}{t} \Big].
	\ee
	Notice also that $\Delta W_{T_1} = W_{T_1}$, $\Delta T_1 = T_1$ and $\Xh_{T_1} := x_0 + b T_1 + W_{T_1}$.
	It follows by Lemma \ref{lemm:AutoDiff_aleaterm} that
	\be \label{eq:Du_toy}
		\partial_x u(0, x_0) 
		~=~ 
		\E \Big[ 
		e^{\beta \Delta T_1} \frac{\Delta W_{T_1}}{\Delta T_1} 
			~\Big( g \big(\Xh_T \big)\1_{\{T \le T_1\}} + \frac{\mu - b}{\beta} \partial_x u \big(T_1,\Xh_{T_1}\big) \1_{\{T_1 < T\}} 
			\Big)
		\Big].~~
	\ee
	Changing the initial condition $(0,x_0)$ to $(T_1, \Xh_{T_1})$, 
	one obtains that, whenever $T_1 < T$,
	$$
		\partial_x u(T_1, \Xh_{T_1}) 
		= 
		\E \Big[ 
		e^{ \beta \Delta T_2} \frac{\Delta W_{T_2}}{\Delta T_2} 
			\Big( g \big(\Xh_T \big)\1_{\{T \le T_2\}} + \frac{\mu - b}{\beta}  \partial_x u \big(T_2,\Xh_{T_2}\big) \1_{\{T_2 < T\}} 
			\Big)
		\Big| T_1, \Xh_{T_1}
		\Big].
	$$
	Plugging the above expression of $\partial_x u(T_1, \Xh_{T_1})$ into the r.h.s. of \eqref{eq:Expr1_toy},
	and using the fact that $T \le T_2$ is equivalent to $N_{T-} \le 1$, and $\P \big[ \{N_{T-} \le 1\} \setminus \{N_T \le 1\} \big] = 0$, 
	it follows that $u(0,x_0) = \E[ \psi_1]$.
	
	\noindent \rmiii Next, changing the initial condition in \eqref{eq:Du_toy} from $(0, x_0)$ to $(T_2, \Xh_{T_2})$ when $T_2 < T$,
	and then plugging the corresponding expression of $\partial_x u(T_2, \Xh_{T_2})$ into $\psi_1$, 
	it follows that $u(0, x_0) = \E[ \psi_2]$.
	Repeating the procedure, we have for all $n \ge 0$,
	$$
		\E \big[ g(X_T) \big] ~=~ u(0, x_0) ~=~ \E [\psi_n ].
	$$
	Finally sending $n \to \infty$, and using  Lemma \ref{lemm:OrderStat_p} together with the dominated convergence theorem,
	it follows that $\E \big[ g(X_T )\big] = \E[\lim_{n \to \infty} \psi_n] = \E [\psi ]$.
	\qed

	\begin{Remark}
		We can also interpret formally the representation $\psi$ in \eqref{eq:psit_vol_constant}
		as the expansion of the diffusion process $X$ around a Brownian motion.
		Let $b=0$ and $ \mu(t,x) \equiv \mu_0$ for some constant $\mu_0 \in \R$,
		so that $\Xh_t = W_t$ and $X_t= \mu_0 t + W_t$.
		Using the fact that $\P(N_T = k) = e^{-\beta T} \frac{(\beta T)^k}{k!}, ~\forall k \ge 0$, 
		it follows formally that
		\be \label{eq:expansion_BM}
			\E \big[ g(X_T) \big]
			=
			\E \left[ \sum_{k = 0}^{\infty} \frac{(\mu_0 T)^k}{k !} g^{(k)}(W_T) \right]
			=
			\E \left[ e^{\beta T}  ~g(W_T) ~\Pi_{k=1}^{N_T} 
				\left( \frac{\mu_0 \Delta W_{k+1}}{\beta \Delta T_{k+1}} \right)
			\right],
		\ee
		where the second equality follows by the fact that 
		$$
			\P[N_T = k] = e^{\beta T} \frac{(\beta T)^k}{k !},
			~\mbox{and}~~
			\E \Big[ g^{(k)} \Big(\sum_{i=0}^k \Delta W_{i+1} \Big) \Big] 
			= 
			\E \Big[ g \Big(\sum_{i=0}^k \Delta W_{i+1} \Big) \Pi_{i=1}^k \frac{\Delta W_{i+1}}{\Delta T_{i+1}} \Big].
		$$
		In particular, the r.h.s. of \eqref{eq:expansion_BM} is exactly $\E \big[\psi \big]$ defined by \eqref{eq:psit_vol_constant} in this  case.
	\end{Remark}

	To conclude this part, we notice that $\psi$ in \eqref{eq:psit_vol_constant} is integrable, but has an infinite variance in general. In the next subsection, we exploit the arbitrariness of the constant $b$ which is involved in the definition of $\psi$. More precisely, we shall choose different constants $b$ at each time $T_k$, in an adaptive way. This will lead to the estimator $\psih$ in Theorems \ref{theo:renorm_constant_vol} and \ref{theo:renorm_general_vol}.

\subsection{A regime switching diffusion representation}

	For $d \ge 1$, $T > 0$, let $(\mu, \sigma): [0,T] \x \R^d \to \R^d \x \M^d$ be bounded continuous functions
	satisfying
	\be \label{eq:LipsCond}
		\big| \mu(t,x) - \mu(t,y) \big|
		+
		\big| \sigma(t,x) - \sigma(t,y) \big|
		~\le~
		L |x-y|;
		~~(t,x,y) \in [0,T] \x \R^d \x \R^d,
	\ee
	for some constant $L > 0$.
	We start by considering a linear parabolic PDE
	\be \label{eq:PDE}
		\partial_t u ~+~ \mu \cdot D u ~+~ a : D^2 u
		~=~ 0,
		&\mbox{on}&
		[0,T) \x \R^d,
	\ee
	with terminal condition $u(T,x) = g(x)$,
	where $a(\cdot) := \frac12\sigma \sigma^{\top}(\cdot)$, $A:B := \mathrm{Tr}(AB^{\top})$ for any two $d\x d$ dimensional matrices $A, B \in \M^d$, and $D,D^2$ denote the gradient and Hessian operators with respect to the space variable $x$.
	Next, let us consider the diffusion process $(X^{0,x_0}_s)_{s \in [0, T]}$
	defined as unique strong solution of the SDE
	\be  \label{eq:SDE1}
		X_0 = x_0,
		&\mbox{and}&
		dX_s
		\;=\;
		\mu \big(s, X_s \big) ~ds
		+
		\sigma \big(s, X_s\big) ~dW_s,
		~~s \in [0,T].
	\ee
	When PDE \eqref{eq:PDE} admits a classical solution in $C^{1,3}_b([0,T] \x \R^d)$,
	i.e. the collection of all functions $\phi(t,x)$ such that $\phi$, $\partial_t \phi$, $D \phi$, $D^2 \phi$ and $D^3 \phi$ all exit and bounded continuous,
	it follows by Feynmann-Kac formula that $V_0 := \E[ g(X^{0,x_0}_T)] = u(0,x_0) $.

	\begin{Remark}
		For technical reason, we will assume that $u \in C^{1,3}_b([0,T] \x \R^d)$ rather than in $C^{1,2}_b([0,T] \x \R^d)$.
		But by approximating the coefficient $\mu, \sigma$ with smooth functions, 
		one can relax this regularity condition in more concrete context. 
	\end{Remark}

	Recall that for $\beta > 0$,
    $(\tau_i)_{i > 0}$ is a sequence of i.i.d. $\Ec(\beta)$-exponential random variables,
    which is independent of the Brownian motion $W$.
    We define
    \b*
        T_k ~:=~ \Big( \sum_{i=1}^k \tau_i \Big) \wedge T,~k\ge 0,
        &\mbox{and}&
        N_t ~:=~ \max \big\{ k : T_k < t \big \}.
    \e*
    Then $(N_t)_{0 \le t \le T}$ is a Poisson process with intensity $\beta$ and arrival times $(T_k)_{k> 0}$,
    and $T_0  = 0$.
    We also introduce, for all $k > 0$, $ \Delta W^k_t ~:=~ W_{(T_{k-1} +t ) \wedge T_k} - W_{T_{k-1}}$.
    It is clear that the sequence of processes $(\Delta W^k_{\cdot})_{k > 0}$ are mutually independent.

	Let $(\muh, \sigmah) : (s,y,t,x) \in  [0,T] \x \R^d  \x [0,T] \x \R^d \longrightarrow \R^d \x \M^d$
	be uniformly bounded, and continuous in $t$, Lipschitz in $x$,
	we define $\Xh$ by
    \be \label{eq:Xh}
        \Xh_0 := x_0
        &\mbox{and}&
        d\Xh_t ~=~ \muh(\Theta_t, t, \Xh_t) dt ~+~ \sigmah(\Theta_t, t, \Xh_t) dW_t,
    \ee
    with $\Theta_t := (T_{N_t}, \Xh_{T_{N_t}})$.
    In other words, the process $\Xh$ is defined recursively by, $\Xh_0 = x_0$ and for all $k \ge 0$,
    \b*
        \Xh_{T_{k+1}}
        =
        \Xh_{T_k}
        +
        \int_{T_k}^{T_{k+1}}
        \muh \big(T_k, X_{T_k}, s , \Xh_s \big) ds
        +
        \int_{T_k}^{T_{k+1}}
        \sigmah \big(T_k, X_{T_k}, s, \Xh_s \big) dW_s.
    \e*
    
	\begin{Example} \label{example:PoisMa}
		\rmi Let $(\muh, \sigmah) (s, y,t,x) = (\mu, \sigma)(s,y)$, then $\Xh$ is defined as a Euler scheme as in \eqref{eq:SDE_Euler_sol}, i.e. $\Xh_0 = x_0$, and
		$$
			\Xh_{T_{k+1}} ~=~ \Xh_{T_k} + \mu(T_k, \Xh_{T_k}) \Delta T_{k+1} + \sigma(T_k, \Xh_{T_k}) \Delta W_{T_{k+1}}.
		$$
		
		\vspace{1mm}
		
		\noindent \rmii When $\muh(\cdot) \equiv 0$ and 
		$\sigmah(s, y, t,x) = \sigma(s,y) + \partial_x \sigma(s,y) (x - y)$,
		then SDE \eqref{eq:Xh} turns to be a linear SDE, whose solution is given explicitly in \eqref{eq:linSDE_sol1}.
	\end{Example}

	We first formulate an assumption on the existence of automatic differentiation  weights associated to SDE \eqref{eq:Xh}.
 	Let $\theta \in [0,T) \x \R^d$ and $(t,x) \in [0,T] \x \R^d$,  the process $(\Xt^{t,x,\theta}_s)_{s \in [t,T]}$ is defined by SDE
	\be \label{eq:EDS_Xt}
	    	\Xt^{t,x,\theta}_t := x,
		~~~
		d \Xt^{t,x,\theta}_s
		~=~
		\muh \big(\theta, s , \Xt^{t,x,\theta}_s \big) ds
		~+~
		\sigmah \big(\theta, s, \Xt^{t,x,\theta}_s \big) dW_s,
	\ee
	\begin{Assumption} \label{assume:MalliavinWeight}
		There is a pair of measurable functions $\big(\Wch^1_{\theta}(\cdot), \Wch^2_\theta(\cdot) \big)$,
		called automatic differentiation weights, taking values in $\R^d \x \M^d$,
		such that, for all $ \theta\in [0,T) \x \R^d$, $(t,x) \in [0,T) \x \R^d$, $s > t$,
		one has $\big(\Wch^i_\theta \big(t, x,  s- t, (W_r -W_t)_{r \in [t,s]}), ~i=1,2 \big)$ are both integrable.
		Moreover, for all bounded continuous function $\phi: \R^d \to \R$,
		\b*
			D^i \E \big[ \phi \big(\Xt^{t,x,\theta}_s \big) \big]
			~=~
			\E \Big[
			\phi \big( \Xt^{t,x,\theta}_s \big)
			~\Wch^i_\theta \big(t, x,  s-t ,  (W_r -W_t)_{r \in [t,s]} \big)
			\Big],
			~~i=1,2,
		\e*
		where $D,D^2$ denote the gradient and Hessian operators with respect to the variable $x$.
	\end{Assumption}

	Let $a (\cdot) := \frac{1}{2} \sigma \sigma^{\top} (\cdot)$ and $\ah (\cdot) := \frac{1}{2} \sigmah \sigmah^{\top} (\cdot)$,
	we denote
	$$
		\Thetah_0 = (t, x)
		~~~\mbox{and then}~~
		\Thetah_k = ( T^t_k, \Xh_{T_k}),
		~~\mbox{for all}~ k > 0.
	$$
	and then for $k > 0$,
	$$
		\Delta f_k
		~~:=~~
		(\mu, a) \big( T_k, \Xh_{T_k} \big)
		-
		(\muh, \ah) \big(\Thetah_{k-1}, T_k, \Xh_{T_k} \big)
		~\in~ \R^d \x \M^d,
	$$
	and for $k \ge 0$,
	$$
		\Wch_k
		~:=~
		\big( \Wch^1_{\Thetah_k}, \Wch^2_{\Thetah_k} \big)
		\big( T_k, \Xh_{T_k}, T_{k+1}, \Delta W^{k+1}_{\cdot} \big)
		~\in~ \R^d \x \M^d,
	$$
    with the weight functions $ \big( \Wch^1_{\theta}(\cdot), \Wch^2_{\theta}(\cdot) \big)$ given in Assumption \ref{assume:MalliavinWeight}.
    We then define
    \be \label{eq:def_general_psih}
        \psih
        ~:=~
        e^{\beta T}
        ~\Big( g(\Xh_T) - g(\Xh_{T_{N_T}}) \1_{\{N_T> 0 \}} \Big)
        ~\beta^{- N_T}~\prod_{k=1}^{N_T}
        \big( \Delta f_k \bullet \Wch_k \big),
    \ee
    where $(p,P) \bullet (q, Q) := p \cdot q + P : Q$ for all $p, q \in \R^d, P, Q \in \M^d$.
    Here we use the convention $\Pi_{k=1}^0 = 1$.
	Finally, for all $n \ge 1$, we also introduce
	\be \label{eq:def_psih_n}
		\psih_n
		&\!\!=&\!\!
		e^{\beta T_{n+1}}
		\prod_{k=1}^{N_T \wedge n} \big(\beta^{-1} \Delta f_k \bullet \Wch_k  \big)
		~\Big[
		\Big( g \big(\Xh_T \big) - g \big(\Xh_{T_{N_T}} \big) \1_{ \{N_T > 0 \}}\Big)
		\1_{\{ N_T \le n \}} \nonumber \\
		&&~~~~~~~~~~~~~~~~~~~~~~~ 
		+ \beta^{-1} \Big(\Delta f_{n+1} \bullet \big( Du, D^2 u \big) \big(T_{n+1}, \Xh_{T_{n+1}} \big) \Big)
		\1_{\{N_T > n\}}
		\Big]. 
	\ee

	\begin{Assumption} \label{assum:ui}
		\rmi The sequence $(\psi_n)_{n \ge 0}$ is uniformly integrable.
		
		\noindent \rmii Let $(e_i)_{i=1,\cdots,d}$ denote the canonical basis of $\R^d$.
		There is some $\eps_0>0$, 
		such that for all $(t,x) \in [0,T)\x\R^d$ and $\theta \in [0,T)\x \R^d$, $n \ge 0$ and $i=1, \cdots,d$, 
		all the following random vectors is integrable:
		$$ 
			\Wch^1_{\theta}(t,x, \tau_1 \wedge (T-t), (W_r - W_t)_{r \in [t, (t+\tau_1) \wedge T]}),
		$$
		$$
			\sup_{\eps \in (0, \eps_0]} \frac{1}{\eps} 
			\Big[ \Wch^1_{\theta} \big( t, x+ \eps e_i, \tau_1 \wedge (T-t), (W_\cdot - W_t) \big) 
				- 
				\Wch^1_{\theta} \big(t, x, \tau_1 \wedge (T-t), (W_\cdot - W_t) \big) \Big]
		$$
		and
		$$  
			\Delta f_{n+1} \bullet \big( Du, D^2 u \big) \big(T_{n+1}, \Xh_{T_{n+1}} \big) \Wch_n.
		$$
	\end{Assumption}

	\begin{Theorem} \label{theo:muh_sigmah}
		Suppose that the PDE \eqref{eq:PDE} has a classical solution $u \in C_b^{1,3} \big([0,T] \x \R^d\big)$,
		suppose in addition that
		Assumptions \ref{assume:MalliavinWeight} and \ref{assum:ui} hold true.
		Then $\psih$ is integrable and $u(0, x_0) = \E \big[\psih\big]$.
	\end{Theorem}

	\begin{Remark} \label{rem:generalTHM}
        {\rm
         \rmi The condition that $u \in C_b^{1,3} \big([0,T] \x \R^d \big)$ may be relaxed in the concrete applications of Theorem \ref{theo:muh_sigmah}.
        This will be indeed performed in Section \ref{theo:renorm_general_vol} by exploiting the integrability of the automatic differentiation  weights $\big(\Wch^1_{\theta}, \Wch^2_{\theta}\big)$ of Assumption \ref{assume:MalliavinWeight}.

        \vspace{1mm}

        \noindent \rmii By definition, the automatic differentiation weight satisfies $\E \big[ \Wch_k \big] = 0$,
        then $\psih$ in \eqref{eq:def_general_psih} has the same mean than the estimator
        \b*
            e^{\beta T}
            ~ g(\Xh_T)
            ~\beta^{- N_T}~\prod_{k=1}^{N_T}
            \big( \Delta f_k \bullet \Wch_k \big).
        \e*
        However, in practice, the weight function $\Wch_k $ is typically of infinity variance, or even not integrable, in general.
        Indeed, as we will see in the following, 
        $\Wch_k$ is generally of order $\frac{1}{\Delta T_{k+1}} = \frac{1}{T_{k+1} - T_k}$,
        where conditioning on $N_T = n$, $(T_1, \cdots, T_{N_T})$ follows the law of statistic order of uniform distribution on $[0,T]$.
        Then by direct computation, one knows $\E \big[ 1/ \Delta T_{N_T+1} \big] = \infty$.
        In the definition of $\psih$ in \eqref{eq:def_general_psih}, the additional term $- g \big(\Xh_{T_{N_T}} \big) \1_{\{N_T> 0 \}} $ 
        can be seen as a control variate so as to guarantee the integrability of $\psih$.
        
	\noindent \rmiii As a consequence of the integrability problems raised in $\mathrm{(ii)}$, Assumption \ref{assum:ui} is in fact implicitly a restriction on the choice of the coefficients $\muh$ and $\sigmah$,
	and we cannot expect a representation for $u(t,x)$ with arbitrary $\muh$ and $\sigmah$,
	see Section \ref{subsec:Proofs} below.
    }
    \end{Remark}

	\noindent {\bf Proof of Theorem \ref{theo:muh_sigmah}.}
	\rmi Recall that $u \in C_b^{1,3} \big([0,T] \x \R^d\big)$ is a classical solution of PDE \eqref{eq:PDE}.
	Denote $(\muh_{\theta}, \ah_{\theta}) (\cdot) = (\muh, \ah)(\theta, \cdot)$,
	one can rewrite \eqref{eq:PDE} in the following equivalent way:
	\be \label{eq:PDE_replace_coef}
		-\partial_t u - \muh_{\theta} \cdot D  u - \ah_{\theta} : D^2 u
		~-~
		\big( (\mu - \muh_{\theta}) \cdot D u + (a  - \ah_{\theta}) : D^2 u \big)
		~=~ 0~.
	\ee	
	Using Feynmann-Kac formula, it follows that
	\be \label{eq:repres_u_0}
		u(0,x_0) 
		=
		\E \Big[ g \big( \Xt^{0,x_0,\theta}_T \big) + \int_0^T \big( (\mu - \muh_{\theta}) \cdot D u + (a  - \ah_{\theta}) : D^2 u \big) \big(s, \Xt^{0,x_0,\theta}_s \big) ds \Big],~~
	\ee
	where $\Xt^{0,x_0,\theta}$ is defined by \eqref{eq:EDS_Xt},
	which coincides with $\Xh$ in \eqref{eq:Xh} on $[0, T_1]$ whenever $\theta = (0, x_0)$.

	Recall that $T_1 = \tau_1 \wedge T$, where $\tau_1$ is a random variable of density $\beta e^{-\beta s} \1_{\{s \ge 0\}}$ independent of the Brownian motion $W$.
	Fixing $\theta = (0, x_0)$, it follows that
        $$
            u(0,x_0)
            \!=\!
            \E \Big[
                e^{\beta T_1}
                \Big( g \big(\Xh_T \big) \1_{\{ N_T = 0 \}} 
                +
                \beta^{-1} \Delta f_1 \bullet  (Du, D^2u) \big(T_1, \Xh_{T_1} \big) \1_{\{N_T > 0\}}
            \Big)
            \Big] 
            \!=\! \E [\psih_0]. 
        $$

	\noindent \rmii Let us now go back to the expression \eqref{eq:repres_u_0}, and derive an expression for 
	the derivatives $D u(0,x_0)$ and $D^2 u(0, x_0)$.
	First, for $D u(0,x_0)$, we use the integrability condition in Assumption \ref{assum:ui} with Lemma \ref{lemm:AutoDiff_aleaterm},
	and also the fact the $Du(\cdot)$ is continuous,
	it follows that
	\b*
		Du(0, x_0) 
		&=& 
		\E \Big[ g \big( \Xt^{0,x_0,\theta}_T \big) \Wch^1_{\theta}(x_0,T) \\
		&&~~~
		+ \int_0^T \big( (\mu - \muh_{\theta}) \cdot D u + (a  - \ah_{\theta}) : D^2 u \big) \big(s, \Xt^{0,x_0,\theta}_s \big) \Wch^1_{\theta}(x_0, s)
		ds \Big],
	\e*
	where we simplify the notation $\Wch^1_{\theta}(0,x_0, s, (W_r - W_t)_{r \in [0, s]})$ to $\Wch^1_{\theta}(x_0, s)$.
	Then by the independence of $\tau_1$ to the Brownian motion $W$, and setting $\theta=(0, x_0)$, 
	it follows that
	\be \label{eq:expression_Du}
		Du(0,x_0) ~=~ \E \Big[ \psih_0 ~\Wch^1_{(0,x_0)}(0,x_0, T_1, \Delta W^1_{\cdot}) \Big].
	\ee
	
	Next, for $D^2u(0,x_0)$, we use again Lemma \ref{lemm:AutoDiff_aleaterm} together with
	the integrability condition in Assumption \ref{assum:ui} and Lipschitz property of 
	$x \mapsto ((\mu - \muh)\cdot Du + (a - \hat a_{\theta}):D^2 u)(s, x)$, and the continuity of $D^2u(\cdot)$ 
	that
	$$
		D^2 u(0,x_0)
		=
		D^2_{x_0} \E \big[ g\big(\Xt^{0,x_0,\theta}_T \big) \big]
		+\int_0^T \!\!D^2_{x_0} \E \Big[ \big((\mu-\muh_{\theta})\cdot Du + (a-\hat a_{\theta}): D^2u\big)(s, \Xt^{0,x_0,\theta}_s) \Big] ds.
	$$
	Setting $\theta = (0, x_0)$ and using Assumption \ref{assume:MalliavinWeight}, it leads to
	$$
		D^2 u(0,x_0) ~=~
		\E \Big[  \psih_0   ~\Wch^2_{(0,x_0)}(0,x_0, T_1, \Delta W^1_{\cdot}) \Big].
	$$
	Recall that $\E [ \Wch^2_{(0,x_0)}(0,x_0, T, \Delta W_{\cdot}) ] = 0$, we then obtain
	\be \label{eq:expression_D2u}
		D^2 u(0,x_0) ~=~
		\E \Big[ \big( \psih_0 - e^{\beta T} g(x_0)\1_{\{N_T= 0\}} \big)  ~\Wch^2_{(0,x_0)}(0,x_0, T_1, \Delta W^1_{\cdot}) \Big].
	\ee
	
	\noindent \rmiii Changing the initial condition $(0,x_0)$ in \eqref{eq:expression_Du} and \eqref{eq:expression_D2u}  by $(T_1, \Xh_1)$ (remember $\psih_0$ dependent also on the initial condition $(0,x_0)$),
	then plugging the expression of $D^1 u(T_1, \Xh_1)$ and $D^2 u(T_1, \Xh_1)$ into the definition of $\psih_0$ in \eqref{eq:def_psih_n},
	it follows by identifying the term that
	\b*
		u(0,x_0) ~=~ \E[ \psi_1].
	\e*

	\noindent \rmiv Repeating the  arguments 
	by replacing the initial condition $(0,x_0)$ by $(T_{n+1}, \Xh_{n+1})$ in \eqref{eq:expression_Du} and \eqref{eq:expression_D2u} 
	and then plugging the corresponding expression into the definition of $\psi_n$, etc.,
	we obtain that $u(0,x_0) = \E[ \psih_n ]$ for all $ n\ge 0$.
	Then letting $n \longrightarrow \infty$, we obtain
	\b*
		u(0,x_0) ~~=~~ \lim_{n \to \infty} \E \big[  \psih_n \big] 
		~~=~~  \E \big[  \lim_{n \to \infty} \psih_n \big]  ~~=~~ \E \big[ \psih \big],
	\e*
	which concludes the proof.
	\qed

\subsection{Proof of the representation results in Theorems \ref{theo:renorm_constant_vol}, \ref{theo:path_depend} and \ref{theo:renorm_general_vol}.}
\label{subsec:Proofs}

	Using the results in Theorem \ref{theo:muh_sigmah}, 
	we can easily complete the proof of the representation results in 
	Theorems \ref{theo:renorm_constant_vol}, \ref{theo:path_depend} and \ref{theo:renorm_general_vol}.

	\vspace{1mm}

	\noindent {\bf Proof of Theorems \ref{theo:renorm_constant_vol} $\mathrm{(ii)}$ and \ref{theo:renorm_general_vol} $\mathrm{(ii)}$.}
	\rmi In the context of Theorems \ref{theo:renorm_constant_vol} and \ref{theo:renorm_general_vol},
	the increment $\Xh_{T_{k+1}} - \Xh_{T_k}$, conditional on $(T_k, \Xh_{T_k})$, is Gaussian.
	And the estimator $\psih$ corresponds to the estimator in Theorem \ref{theo:muh_sigmah}
	with automatic differentiation weights function
	\be \label{eq:MalliavinW_0}
		\Wch^1_{\theta} \big(\cdot, \delta t, \delta w \big) := (\sigma^{\top}_0)^{-1} \frac{\delta w}{\delta t}
		&\mbox{and}&
		\Wch^2_{\theta} \big(\cdot,\delta t, \delta w \big) := (\sigma^{\top}_0)^{-1} \frac{\delta w \delta w^{\top} - \delta t I_d}{\delta t^2} \sigma_0^{-1}.~~~~
	\ee
	In particular, it is clear that Assumption \ref{assume:MalliavinWeight} holds true with the above choice of  automatic differentiation  weight functions in \eqref{eq:MalliavinW_0}.
	
	\noindent \rmii Besides, the uniform integrability conditions and integrability conditions in Assumption \ref{assum:ui}
	can be easily obtained following the lines in the first part of the proof of Theorems \ref{theo:renorm_constant_vol} and \ref{theo:renorm_constant_vol}, using Lemma \ref{lemm:OrderStat_p}.

	\noindent \rmiii Now, suppose in addition that $\mu$, $\sigma$ and $g$ are bounded smooth functions with bounded continuous derivatives,
	so that $u \in C_b^{1,3}([0,T] \x \R^d)$.
	It follows by Theorem \ref{theo:muh_sigmah} that $V_0 = \E[ \psih]$.

	\vspace{1mm}

	\noindent \rmiv  Finally, when $\mu(\cdot) $ and $\sigma(\cdot) $ satisfy the Lipschitz condition \eqref{eq:Lipschitz_coef2}
    and $g$ is Lipschitz, 
    we can find a sequence of bounded smooth functions $(\mu_{\eps}(\cdot) , \sigma_{\eps}(\cdot) , g_{\eps}(\cdot) )$ which converges locally uniformly to $(\mu(\cdot) , \sigma(\cdot) , g(\cdot) )$ as $\eps \to 0$.
    Let $X^{\eps}$ be the solution of
    \b*
        dX^{\eps}_t &=& \mu_{\eps}(t, X^{\eps}_t) dt ~+~ \sigma_{\eps}(t, X^{\eps}_t) dW_t.
    \e*
    Then by the stability of SDEs together with dominated convergence theorem, it follows that
    \b*
        V^{\eps}_0 ~:=~ \E \big[ g_{\eps}(X^{\eps}_T) \big]
        ~&\longrightarrow&~
        V_0 := \E \big[ g(X_T) \big],
        ~~~~\mbox{as}~ \eps ~\rightarrow~ 0.
    \e*
    Moreover, by Lemma \ref{lemm:OrderStat_p} together with dominated convergence theorem,
    it is easy to prove that $\E[ \psih^{\eps}]  \to \E[ \psih] $ as $\eps \to 0$,
    where $\psih^{\eps}$ denotes the estimator of the algorithm \eqref{eq:def_psih2} associated to the coefficient $(\mu_{\eps}, \sigma_{\eps}, g_{\eps})$.
    We then conclude the proof.
    \qed

	\vspace{2mm}

	\noindent {\bf Proof of Theorem \ref{theo:path_depend} $\mathrm{(ii)}$.}
	For the path-dependent case, it is enough to use the same arguments as in Theorem \ref{theo:renorm_constant_vol},
	together with the PDE system \eqref{eq:PDE_n} in Remark \ref{rem:PDE_n}.
	\qed

\subsection{Proof of Theorem \ref{theo:ExactSimulation2}}
\label{subsec:proof_0drift1d}

    To introduce the algorithm in the context of Theorem \ref{theo:muh_sigmah}, we propose to choose
    \b*
        \muh(\cdot) \equiv 0
        &\mbox{and}&
        \sigmah(s, y, t,x) = \sigma(s,y) + \partial_x \sigma(s,y) (x - y).
    \e*

    Before providing the proof of Theorem \ref{theo:ExactSimulation2}, 
    we first give a lemma which justifies our choice of the automatic differentiation  weight function $\Wcb^2_k$ in \eqref{eq:MWeight2},
    as well as some related estimations.
    Let $c_1, c_2, x \in \R$ be constants such that $ c_1 + c_2 x \neq 0$,
    we denote by $\Xb^{0,x}$ solution of the SDE
    \be \label{eq:sol_SDE_dynamic}
        \Xb_0 = x,
        &&
        d \Xb_t
        ~=~
        \big (c_1 + c_2 \Xb_t \big) dW_t,
    \ee
    whose solution is given explicitly by
    \be \label{eq:sol_SDE}
        \Xb^{0,x}_t
        ~=
        \begin{cases}
            -~
            \frac{c_1}{c_2}
            ~+~ \big( \frac{c_1}{c_2} + x \big) \exp \Big(- \frac{c_2^2}{2} t + c_2 W_t \Big),
            &~\mbox{if}~ c_2 \neq 0, \\
            x + c_1 W_t,
            &~\mbox{if}~ c_2 = 0.
        \end{cases}
    \ee
    Consider also its antithetic variable $\Xt^x_t$ defined by
    \b*
        \Xt^{0,x}_t
        ~=
        \begin{cases}
            -~
            \frac{c_1}{c_2}
            ~+~ 
            \big( \frac{c_1}{c_2} + x \big) 
            \exp \Big(- \frac{c_2^2}{2} t - c_2 W_t \Big),
            &~\mbox{if}~ c_2 \neq 0, \\
            x - c_1 W_t,
            &~\mbox{if}~ c_2 = 0.
        \end{cases}
    \e*

    \begin{Lemma} \label{lemm:MalliavinW_LogNorm}
        Let $x \in \R$, $(c_1, c_2) \in \R^2$ be two constants such that $c_1 + c_2 x \neq 0$,
        $\phi: \R \to \R$ a bounded continuous function. \\
        \noindent \rmi
        Then for all $t \in (0, T]$,
        \be \label{eq:second_derivative_weight}
            \partial^2_{xx} \E \big[ \phi \big(\Xb^{0,x}_t \big) \big]
            &=&
            \E \Big[
                \phi \big(\Xb^{0,x}_t \big)  \frac{1}{(c_1 + c_2 x)^2}
                \Big( - c_2 \frac{W_t}{t} + \frac{W_t^2 - t}{t^2} \Big)
            \Big].
        \ee
        \noindent \rmii Suppose in addition that $\phi (\cdot) \in C^2_b(\R)$.
        Then there is some constant $C$ independent of $(t,x)$ such that,
                for all $(t,x) \in [0,T] \x \R^d$,
        \b*
            &&\E \Big[ \Big( \phi \big(\Xb^{0,x}_t\big) - \phi(x) \Big)^2 \Big( \frac{W_t}{t} \Big)^2 \Big]
            +
            \E \Big[
                \Big( \phi \big(\Xb^{0,x}_t \big) - 2\phi(x) + \phi(\Xt^{0,x}_t) \Big)^2
                \Big( \frac{W_t^2 - t}{t^2} \Big)^2
            \Big] \\
            &\le&
            C (c_1 + c_2 x)^2.
        \e*
    \end{Lemma}
    \proof
    \rmi First, when $c_2 = 0$, it is clear that result is correct (see e.g. Lemma 2.1 of Fahim, Touzi and Warin \cite{fah}).
    Next, when $c_2 \neq 0$, denote $v(x) := \E \big[ \phi \big(\Xb^{0,x}_t \big) \big]$, 
    then with the expression of $\Xb^{0,x}_t$ in \eqref{eq:sol_SDE}, it follows that
    \b*
        v(x)
        &=&
        \int_{\R}
        \phi
            \Big( - \frac{c_1}{c_2} + \Big( \frac{c_1}{c_2} + x \Big) e^{- c_2^2 t/2 + c_2 \sqrt{t} y} \Big)
            \frac{1}{\sqrt{2 \pi}} e^{-y^2/2} dy.
    \e*
    Suppose that $\phi(\cdot) \in C^2_b(\R)$,
    then using integration by parts, it follows that
    \b*
        v'(x)
        &=&
        \int_{\R}
        \phi'
            \Big( - \frac{c_1}{c_2} + \Big( \frac{c_1}{c_2} + x \Big) e^{- c_2^2 t/2 + c_2 \sqrt{t} y} \Big)
            e^{- c_2^2 t/2 + c_2 \sqrt{t} y}
            \frac{1}{\sqrt{2 \pi}} e^{-y^2/2} dy \\
        &=&
        \int_{\R}
        \phi
            \Big( - \frac{c_1}{c_2} + \Big( \frac{c_1}{c_2} + x \Big) e^{- c_2^2 t/2 + c_2 \sqrt{t} y} \Big)
            \frac{1}{c_1+c_2x}  \frac{y}{\sqrt{t}}
            \frac{1}{\sqrt{2 \pi}} e^{-y^2/2} dy \\
        &=&
        \E\Big[ \phi(\Xb^{0,x}_t) \frac{1}{c_1+c_2x} \frac{W_t}{t} \Big].
    \e*
    Similarly, still  using integration by parts, and by direct computation, we obtain
    \b*
        v''(x)
        &=&
        \E \Big[
            \phi \big(\Xb^{0,x}_t \big)  \frac{1}{(c_1 + c_2 x)^2}
            \Big( - c_2 \frac{W_t}{t} + \frac{W_t^2 - t}{t^2} \Big)
        \Big].
    \e*
    When $\phi(\cdot)$ is only a bounded continuous function, one can approximate $\phi(\cdot)$ by a sequence of smooth function $\phi_{\eps}(\cdot)$ which converges to $\phi(\cdot)$ uniformly, and $\phi'_{\eps}$ and $\phi''_{\eps}$ are bounded continuous.
    We then obtain
    $$
        v_{\eps}(x) ~:=~ \E \big[ \phi_{\eps} \big( \Xb^{0,x}_t \big) \big] ~~\to~~ v(x).
    $$
    Moreover, the limit $\lim_{\eps \to 0} v'_{\eps}(x)$, $\lim_{\eps \to 0} v''_{\eps}(x)$ exist,
    thus $v''(x)$ also exists and
    \b*
        v''(x) ~=~ \lim_{\eps \to 0} v''_{\eps}(x) ~=~
        \E \Big[
            \phi \big(\Xb^{0,x}_t \big)  \frac{1}{(c_1 + c_2 x)^2}
            \Big( - c_2 \frac{W_t}{t} + \frac{W_t^2 - t}{t^2} \Big)
        \Big].
    \e*

    \noindent \rmii When $c_2 = 0$, the estimation in \rmii of the statement is clear true since $\phi'$ and $\phi''$ are uniformly bounded.
   
    When $c_2 \neq 0$, denote $| \phi'|_0 := \sup_{x} |\phi'(x)|$, we obtain, by direct computation, that
    \b*
        &&
        \E \Big[ \Big( \phi \big(\Xb^{0,x}_t\big) - \phi(x) \Big)^2 \Big( \frac{W_t}{t} \Big)^2 \Big]
        ~~\le~~
        |\phi'|_0 \E \Big[ \big(\Xb^{0,x} - x \big)^2 \frac{W_t^2}{t^2} \Big] \\
        &=&
        |\phi'|_0
        \E \Big[
            \big(c_1 + c_2 x \big)^2
            \Big(\frac{e^{-c_2^2 t/2 + c_2  W_t} - 1}{c_2 W_t - c_2^2 t/2}\Big)^2
            \frac{W_t^2 (c_2 W_t - c_2^2 t/2)^2}{t^2}
        \Big],
    \e*
    which is clearly uniformly bounded by $C (c_1 + c_2 x)^2$ for some constant $C$ independent of $(t,x) \in [0,T] \x \R^d$. 

    Next, denote $\ell(y) := \big( x + \frac{c_1}{c_2} \big) \big ( e^{- c_2^2 t/2 + c_2 y} - 1 \big) $,
    and define $\varphi(y) := \phi( x + \ell(y))$.
    Then
    \be \label{eq:varphi_2}
        \varphi''(y) = \phi''(x + \ell(y) ) (c_2 + c_1 x)^2 e^{- c_2^2 t + 2 c_2 y}
            + \phi'(x +\ell(y)) (c_2+ c_1 x) c_2 e^{- c_2^2 t/2 + c_2 y}.~
    \ee
    It follows by the definition of $\varphi$ as well as its derivative, together with direct computation,  that
    \b*
        &&
        \E \Big[
            \Big( \phi \big(\Xb^{0,x}_t \big) - 2\phi(x) + \phi(\Xt^{0,x}_t) \Big)^2
            \Big( \frac{W_t^2 - t}{t^2} \Big)^2
        \Big] \\
        &=&
        \E \Big[
            \big(\varphi( W_t) + \varphi( - W_t) - 2 \varphi(0) \big)^2
            \Big( \frac{W_t^2 - t}{t^2} \Big)^2
        \Big]
        ~+~
        \E \Big[
            2 \big( \varphi(0) - \phi(x) \big)^2
            \Big( \frac{W_t^2 - t}{t^2} \Big)^2
        \Big]\\
        &\le&
        \E \Big[
            \Big( \frac{W_t^2 (W_t^2 - t)}{t^2} \Big)^2
           \sup_{|z| \le |W_t|} \varphi''(z)
        \Big] \\
        &&
        ~+~
        \E \Big[
            2 \Big( \phi \Big( x+ \frac{c_1 + c_2 x}{c_2} \big(e^{-c^2_2 t/2 } - 1 \big) \Big) - \phi(x) \Big)^2
            \Big( \frac{W_t^2 - t}{t^2} \Big)^2
        \Big],
    \e*
    which is also uniformly bounded by $ C (c_1 + c_2 x)^2$  for some constant $C > 0$,
    \qed

    \vspace{2mm}

    \noindent {\bf Proof of Theorem \ref{theo:ExactSimulation2}.}
    \rmi Let us first prove that $\E \big[ \psib^2 \big] < \infty$ for $\psib$ defined by \eqref{eq:psib2}.
    First, we notice that $\Wch^-_k = \Wch_k^2$ for all $k = 1, \cdots, N_T-1$, $g \in C^2_b(\R)$,
    and with the choice of $c_1^k$ and $c_2^k$ in \eqref{eq:def_c1c2k}, one has
    $ c_1^k + c_2^k \Xh_{T_k} = \sigma(T_k, \Xh_{T_k})$, which is uniformly bounded.
    By considering the conditional expectation over $(\Xh_{T_{N_T}}, \Delta T_{N_T+1})$
    using items \rmii of Lemma \ref{lemm:MalliavinW_LogNorm}, we have 
    $ \E \big[ \big| \psib \big|^2 \big] $ is bounded by
    \b*
        ~~C \E \Big[
        \beta^{-2 N_T} \prod_{k=2}^{N_T}
        \Big\{
            \frac{a(T_k, \Xh_{T_k}) - \tilde a_k}{2 a(T_k, \Xh_{T_k})}
            \Big(- \partial_x \sigma(T_k, \Xh_{T_k}) \frac{\Delta W_{T_k}}{\Delta T_k}
            + \frac{\Delta W_{T_k}^2 -\Delta T_k}{\Delta T_k^2} \Big)
        \Big \} ^2
        ~\Big],
    \e*
    for some constant $C$.
    Further, by denoting $\Delta \Xh_{T_k} := \Xh_{T_k} - \Xh_{T_{k-1}}$, one has
    \b*
        \big| a(T_k, \Xh_{T_k}) - \tilde a_k \big|
        &\le&
        \Big( |\sigma|_0 + \big |\partial_x \sigma(T_{k-1}, \Xh_{T_{k-1}})\Delta \Xh_{T_k} \big| /2 \Big) \\
        &&~~~~~~~~~~~~~~~~~~~~~
        \Big( |\partial_t \sigma|_0 \Delta T_k + \big| \partial^2_{xx} \sigma \big|_0 \big( \Delta \Xh_{T_k} \big)^2 \Big),
    \e*
    where $|\sigma|_0 := \sup_{t,x} |\sigma(t,x)|$.
    Notice that $\sigma \ge \eps > 0$, $\sigma$ and $\partial_x \sigma$ are uniformly bounded,
    then to prove that $\psib$ is of finite variance, 
    it is enough to prove that, for some $C > 0$ large enough, the expectation of
    \be \label{eq:Product1}
        \prod_{k=2}^{N_T}
        \Big[
            C \Big( C +  \big |\partial_x \sigma(T_{k-1}, \Xh_{T_{k-1}})\Delta \Xh_{T_k} \big| \Big) ^2
            \Big( C + \frac{\Delta \Xh_{T_k}^2}{\Delta T_k} \Big)^2
            \Big( C \big| \Delta W_{T_k} \big| + \frac{\Delta W_{T_k}^2}{\Delta T_k} + 1 \Big)^2
        \Big]~~
    \ee
    is finite.
    Similarly to the computation in item \rmii of Lemma \ref{lemm:MalliavinW_LogNorm}, we have
    \b*
        && \Delta \Xh_{T_k}
        ~~=~~
        \Xh_{T_k} - \Xh_{T_{k-1}} \\
        &=&
        \sigma(T_{k-1}, \Xh_{T_{k-1}})
        \frac{\exp \big( - \partial_x \sigma(T_{k-1}, \Xh_{T_{k-1}})^2 \Delta T_k /2
            + \partial_x \sigma(T_{k-1}, \Xh_{T_{k-1}}) \Delta W_{T_k} \big)
            - 1
        }
        {\partial_x \sigma(T_{k-1}, \Xh_{k-1}) }.
    \e*
    Notice again that $\sigma(\cdot)$ and $\partial_x \sigma(\cdot)$ are uniformly bounded,
    it follows that
    \b*
        &
        \E \Big\{
            \Big[ \Big( C +  \big |\partial_x \sigma(T_{k-1}, \Xh_{T_{k-1}})\Delta \Xh_{T_k} \big| \Big)
            \Big( C + \frac{\Delta \Xh_{T_k}^2}{\Delta T_k} \Big)
            \Big( C \big| \Delta W_{T_k} \big| + \frac{\Delta W_{T_k}^2}{\Delta T_k} + 1 \Big)
        \Big]^2&
        \\
        &~~~~~~~~~~~~~~~~~~~~~~~~~~~~~~~~~~~~~~~~~~~~~~~~~~~~~~~~~~~~~~~~~~
        \Big|~
        ~\Xh_{T_{k-1}}, T_{k-1}, \Delta T_k
        ~\Big \}
        &\le~
        C',
    \e*
    for some constant $C' > 0$ independent of $\Xh_{T_{k-1}}, T_{k-1}, \Delta T_k$.
    Then the variance of \eqref{eq:Product1} is 
    bounded by $C \E \big[ (C')^{N_T} \big] < \infty$ and hence $\psib$ in \eqref{eq:psib2} is of finite variance.

    \vspace{1mm}

    \noindent \rmii Let us now consider the estimator $\psih$.
    By the same computation, we obtain that
    \b*
        \E \big[~ \psih ~\big|~ N_T, \Delta T_1, \cdots, \Delta T_{N_T+1} \big]
        & \le &
        C^{N_T} \frac{1}{\sqrt{\Delta T_{N_T+1}}},
        ~~\mbox{for some}~
        C > 0,
    \e*
    where the r.h.s. is integrable but of infinite variance (see Lemma \ref{lemm:OrderStat_p}).
    Similarly, it is easy to check the uniform integrability condition in Assumption \ref{assum:ui} for $\psih$ in \eqref{eq:psih2}.

    \vspace{1mm}

    \noindent \rmiii Finally, using Lemma \ref{lemm:MalliavinW_LogNorm} \rmi, it follows that Assumption \ref{assume:MalliavinWeight} holds true.
    Moreover, with the regularity condition on $\sigma(t,x)$ and $g$ in Assumption \ref{assum:1d_driftless}, we know $u \in C_b^{1,3}(\R)$.
    We then conclude the proof of $u(0,x_0) = \E [ \psih] = \E[ \psib]$
    by Theorem  \ref{theo:muh_sigmah}.
    \qed

\appendix

\section{Appendix}

	We first provide an estimation on the order statistics of uniform distribution on $[0,1]$,
	which induces an estimation on a functional of the arrival times $(T_k)_{k > 0}$ of the Poisson process.
	We next provide a technical result on the automatic differentiation function related to a SDE.

    \begin{Lemma} \label{lemm:order_stat}
        Let $p \in (0,1)$, $(U_k)_{k = 1, \cdots, m}$ be a sequence of i.i.d. random variable of uniform distribution on $[0,1]$,
        and $(U_{(1)} \le U_{(2)} \le \cdots \le U_{(m)}$ be the associated order statistics.
        Then
        \b*
            G_{m,p} ~~:=~~
            \E \left [ \Big( \frac{1}{U_{(1)}} \frac{1}{U_{(2)}- U_{(1)}}
            \cdots
            \frac{1}{U_{(m)}- U_{(m-1)}} \Big)^p  \right]
            &\le&
            m! ~\frac{1}{(1-p)^m}.
        \e*
    \end{Lemma}
    \proof First, we notice that for any $x \in (0,1)$,
    \be \label{eq:integral}
        \int_x^1 \big( \frac{1}{u - x} \big)^p du
        &=&
        \frac{(1-x)^{1-p}}{1-p}
        ~~\le~~
        \frac{1}{1-p}.
    \ee
    Then, since the density of the order statistics $(U_{(1)}, \cdots, U_{(m)})$ is provided by
    \b*
        f(u_1, \cdots, u_m) &:=& m! ~ \1_{\{0 < u_1 < u_2 < \cdots < u_m < 1\}},
    \e*
    it follows by direct computation that
    \b*
        G_m =
        m! ~\int_0^1 \int_{u_1}^1 \cdots \int_{u_{m-1}}^1 \Big( \frac{1}{u_1} \frac{1}{u_2 - u_1} \cdots \frac{1}{u_m - u_{m-1}} \Big)^p ~du_1 \cdots d u_m
        \le m!~ \frac{1}{(1-p)^m},
    \e*
    where the last inequality follows from \eqref{eq:integral}. \qed

    \vspace{3mm}

	Let $N = (N_s)_{s \ge 0}$ be a Poisson process with arrival times $(T_k)_{k > 0}$, denote $\Delta T_{k+1} := T_{k+1} - T_k$.
	Let $0 = t_0 < t_1 < \cdots < t_n = T < \infty$ be a discrete time grid, we define further
	\b*
		\Tt_k~ :=~ \min(T_k, t_i),
		&\mbox{whenever}&
		T_{k-1} \in [t_{i-1}, t_i)
		~\mbox{for some}~
		i = 1, \cdots, n,
		\e*
	and
	\b*
    	        \Delta \Tt_k ~:=~ \Tt_k - \Tt_{k-1}
	        &\mbox{for every}&
	        k = 2, \cdots, N_T +1.
	\e*

    \begin{Lemma} \label{lemm:OrderStat_p}
        Let $p \in (0,1)$, then for every constant $C > 0$, one has
        \b*
            \E \Big[~
            	\prod_{k=1}^{N_T} 
		\frac{C}{(\Delta \Tt_{k+1})^p} 
           ~ \Big]
            &<& \infty.
        \e*
    \end{Lemma}
    \proof \rmi    Notice that we can always add points into the time grid $0 = t_0 < t_1 < \cdots < t_n = T$,
    which makes $\Delta \Tt_{i+1}$ smaller.
    Therefore, one can suppose without loss of generality that $t_k - t_{k-1} < (2 \beta C)^{p-1}$ for every $k = 1, \cdots, n$.

    \vspace{1mm}

    \noindent \rmii For every $k = 1, \cdots, n$, we denote $N^k := \# \{i ~: T_i \in [t_{k-1}, t_k)\}$,
    and $\Tt^k_i := \Tt_{k_i}$ with $k_i := \sum_{j < k} N^j + i$ for $i = 1, \cdots, N^k+1$,
    and $\Delta \Tt^k_i := \Tt^k_i - \Tt^k_{i-1}$.
    By the memoryless property of the exponential distribution,
    it is clear that $\big(\Delta \Tt^1_i, ~i = 2, \cdots N^1+1 \big), \cdots, \big(\Delta \Tt^n_i, ~i = 2, \cdots N^n+1 \big)$ are mutually independent.
    Moreover, we have
    \be \label{eq:Prod_interm}
        \prod_{i=1}^{N_T} ~\frac{C}{ (\Delta \Tt_{i+1})^p}
        &=&
        \prod_{k=1}^{n} ~\left( \prod_{i=1}^{N^k} ~\frac{C}{(\Delta \Tt^k_{i+1})^p} \right).
    \ee
    Next, the law of $\big(T^k_i, ~i = 1, \cdots, N^k\big)$ conditioning on $N^k = m$ is the law of order statistics of uniform distribution on $[t_{k-1}, t_{k}]$.
    Then it follows by Lemma \ref{lemm:order_stat} that for every $k = 1, \cdots, n$,
    \b*
        \E\left[~ \prod_{i=1}^{N^k} ~ \frac{C}{(\Delta \Tt^k_{i+1})^p} ~\right]
        &\le&
        e^{-\beta (t_{k}- t_{k-1})}
        \sum_{m=0}^{\infty} \frac{ (\beta (t_{k}- t_{k-1}))^m }{m!}
            m! 2^m \left( \frac{ C}{(t_{k}- t_{k-1})^p} \right)^m \\
        &=&
        e^{-\beta (t_{k}- t_{k-1})}
        \sum_{m=0}^{\infty} \left(2 \beta C (t_{k} - t_{k-1})^{1-p} \right)^m
        ~<~ \infty,
    \e*
    where the last inequality follows by the fact $t_{k} - t_{k-1} <  (2 \beta C)^{p-1}$.
    We then conclude the proof by \eqref{eq:Prod_interm}.
    \qed

	\vspace{3mm}
	
	Let $X^x$ be the solution of SDE
	\b*
		X^x_0 = x,
		~~~
		dX^x_t ~=~ \mu(t, X_t^x) dt ~+~ \sigma(t, X^x_t) dW_t,
	\e*
	where $(\mu, \sigma): [0,T] \x \R^d \to \R^d \x \M^d$ is continuous and in addition Lipschitz continuous in $x$.
	\begin{Lemma} \label{lemm:AutoDiff_aleaterm}
		Suppose that for all $t \in [0,T]$ and bounded continuous function $\phi: [0,T] \x \R^d \to \R$,
		the derivatives 
		$\big( \partial_{x_i}  \E \big[ \phi(t, X^x_t) \big], \partial^2_{x_i,x_j}  \E \big[ \phi(t, X^x_t) \big] \big)_{i,j =1,\cdots,d}$ exist;
		and there is some measurable $\R^d$-valued function $\Wch^1 \big(x, t, (W_s)_{s\in [0,t]} \big)$
		such that for 
		$$
			\partial_{x_i} 
			\E \big[ \phi(t, X^x_t) \big]  ~=~ \E \big[ \phi(t, X^x_t) \Wch^1_i(x, t, (W_s)_{s\in [0,t]}) \big],
			~~~i=1,\cdots, d.
		$$
		Let $F(dt)$ be some probability measure on $[0,T]$.
		
		\noindent \rmi Suppose that
		for each $x \in \R^d$, $i= 1, \cdots, d$,
		$$
			\int_0^T \E\Big[  \big|  \Wch^1_i(x, t, (W_s)_{s\in [0,t]}) \big|  \Big] F(dt)
			~<~
			\infty,
		$$
		and the continuous function $\phi: [0,T] \x \R^d \to \R$ (which may be unbounded) satisfies
		$$
			\int_0^T \E\Big[  \big| \phi(t,X^x_t)  \Wch^1_i(x, t, (W_s)_{s\in [0,t]}) \big|  \Big] F(dt)
			~<~
			\infty.
		$$
		Then
		\be \label{eq:partial1}
			\partial_{x_i} ~ \int_0^T  \E \big[\phi(t, X^x_t)  \big]  F(dt)
			~=~
			\int_0^T\E \Big[  \phi(t, X^x_t) \Wch^1_i (x, t, (W_s)_{s\in [0,t]})  \Big] F(dt) .
		\ee
		
		\noindent \rmii Suppose in addition that $\phi(t,x)$
		is bounded continuous function and Lipschitz in $x$,
		and for each $x \in \R^d$, $i,j = 1, \cdots d$,
		$$
			\int_0^T \sqrt{ \E\Big[  \big|  \Wch^1_i(x, t, (W_s)_{s\in [0,t]}) \big|^2  \Big]} F(dt)
			~<~
			\infty,
		$$	
		and
		\be \label{eq:partial_W1}
			\int_0^T \sup_{\eps \in [0, \eps_0]} \Big| \frac{1}{\eps} 
				\E\Big[ 
					\Wch^1_i (x+ \eps e_j, t, (W_s)_{s\in [0,t]}) - \Wch^1_i (x, t, (W_s)_{s\in [0,t]})
				\Big]
			\Big|
			F(dt)
			~<~
			\infty,
		\ee
		for some $\eps_0 > 0$,
		where $(e_j)_{j=1,\cdots,d}$ denotes the canonical basis of $\R^d$.
		Then
		\be \label{eq:partial2}
			\partial^2_{x_i x_j} ~\int_0^T \E \Big[  \phi(t, X^x_t)\Big]  F(dt) 
			~=~ 
			\int_0^T \partial^2_{x_i x_j} \E \big[ \phi(t, X^x_t)\big] F(dt),
		\ee
		where, in particular, the partial derivative at the l.h.s. and the integration at the r.h.s. are well defined.
	\end{Lemma}
	\proof
	\rmi First, let us notice that $(t,x) \mapsto (\mu, \sigma)(t,x)$ is Lipschitz in $x$, 
	then by standard analysis (see e.g. Chapter 7.8 of \cite{GrahamTalay}), there is some constant $C$ independent of $\eps > 0$ and $i=1, \cdots d$, such that
	\be \label{eq:moment_var1}
		\E \Big[ \Big| \frac{X^{x+\eps e_i} - X^x}{\eps} \Big|^2 \Big] 
		~\le~
		C (1 + e^{Ct}).
	\ee
	
	\noindent \rmii Suppose that $\phi(t,x)$ is bounded continuous and Lipschitz in $x$.
	It follows that
	\b*
		&&
		\lim_{\eps \to 0} \int_0^T \frac{1}{\eps} \E \Big[  \big( \phi(t, X^{x + \eps e_i}_t) - \phi(t, X^x_t) \big)\Big] F(dt) \\
		&=&
		\int_0^T \lim_{\eps \to 0} \frac{1}{\eps} \E \Big[  \big( \phi(t, X^{x + \eps e_i}_t) - \phi(t, X^x_t) \big)\Big] F(dt) \\
		&=& 
		\E \Big[ \int_0^T \phi(t, X^x_t) \Wch^1_i (x, t, (W_s)_{s\in [0,t]}) F(dt)  \Big],
	\e*
	where the first equality follows by the Lipschitz property of $x \mapsto \phi(t,x)$ and \eqref{eq:moment_var1}.
	We hence proved \eqref{eq:partial1} when $x \mapsto \phi(t,x)$ is Lipschitz.

	\vspace{1mm}

	\noindent \rmiii When $\phi$ is only continuous, it is enough to approximate it by a sequence
	$(\phi_n)_{n\ge1}$ which are all bounded, and Lipschitz in $x$.
	Then by the integrability of $ \Wch^1_i(x, t, (W_s)_{s\in [0,t]})$ as well as that of 
	$\phi(t, X^x_t) \Wch^1_i(x, t, (W_s)_{s\in [0,t]})$ under $\P(d\om) \x F(dt)$, 
	it follows that \eqref{eq:partial1} holds true for continuous function $\phi$.
	
	\vspace{1mm}

	\noindent \rmiv 	To prove \eqref{eq:partial2}, let us use \eqref{eq:partial1} and obtain that
	\be \label{eq:lim_int_2}
		&&
		\lim_{\eps \to 0} \frac{1}{\eps} \Big[
			\partial_{x_j} \int_0^T \E \big[ \phi(t, X^{x+\eps e_i}_t) \big]  F(dt)
			-
			\partial_{x_j} \int_0^T \E \big[ \phi(t, X^{x}_t) \big]  F(dt)
		\Big]  \nonumber \\
		&=&
		\lim_{\eps \to 0}  \int_0^T \frac{1}{\eps}
			\E \Big[ 
				\phi(t, X^{x+\eps e_i}_t) \Wch^1_j(x+\eps e_i, t, \cdot)
				- \phi(t, X^x_t) \Wch^1_j(x, t,\cdot)
			\Big] F(dt) \nonumber \\
		&=&
		 \int_0^T \lim_{\eps \to 0}  \frac{1}{\eps}
			\E \Big[ 
				\phi(t, X^{x+\eps e_i}_t) \Wch^1_j(x+\eps e_i, t, \cdot)
				- \phi(t, X^x_t) \Wch^1_j(x, t,\cdot)
			\Big] F(dt).
	\ee
	where the first equality follows by the Lipschitz property of $x \mapsto \phi(t,x)$ and the estimation \eqref{eq:moment_var1} together with  
	\eqref{eq:partial_W1},
	and in particular, the integrable in the last term of \eqref{eq:lim_int_2} is well defined,
	and hence the limit of the first term of \eqref{eq:lim_int_2}  exists.
	\qed

\end{document}